\documentclass[10pt]{article}
\font\smallit=cmti10
\font\smalltt=cmtt10

\usepackage{amssymb,latexsym,amsmath,epsfig,amsthm,float,graphicx,MnSymbol,pgfplots,tikz, url,wasysym}

\usepackage{comment}

\usepackage[utf8]{inputenc}
\pgfplotsset{compat=1.17}

\makeatletter

\renewcommand\section{\@startsection {section}{1}{\z@}
{-30pt \@plus -1ex \@minus -.2ex}
{2.3ex \@plus.2ex}
{\normalfont\normalsize\bfseries\boldmath}}

\renewcommand\subsection{\@startsection{subsection}{2}{\z@}
{-3.25ex\@plus -1ex \@minus -.2ex}
{1.5ex \@plus .2ex}
{\normalfont\normalsize\bfseries\boldmath}}

\renewcommand{\@seccntformat}[1]{\csname the#1\endcsname. }

\makeatother

\newtheorem{theorem}{Theorem}
\newtheorem{lemma}{Lemma}
\newtheorem{proposition}{Proposition}

\newtheorem*{nonum-theorem}{Theorem}

\theoremstyle{definition}
\newtheorem{definition}{Definition}
\newtheorem{conjecture}{Conjecture}
\newtheorem{remark}{Remark}

\newtheorem*{nonum-remark}{Remark}

\makeatletter
\newtheorem*{rep@theorem}{\rep@title}
\newcommand{\newreptheorem}[2]{
\newenvironment{rep#1}[1]{
\def\rep@title{#2 \ref{##1}}
\begin{rep@theorem}}
{\end{rep@theorem}}}
\makeatother

\makeatletter
\newtheorem*{rep@conjecture}{\rep@title}
\newcommand{\newrepconjecture}[2]{
\newenvironment{rep#1}[1]{
\def\rep@title{#2 \ref{##1}}
\begin{rep@conjecture}}
{\end{rep@conjecture}}}
\makeatother

\newreptheorem{theorem}{Theorem}
\newreptheorem{lemma}{Lemma}

\newreptheorem{conjecture}{Conjecture}

\newcommand{\Prob}{\mathbb{P}}

\newcommand{\jstar}{j^{\ast}}

\newcommand\shortunderline[2][3]{\mkern#1mu\underline{\mkern-#1mu#2\mkern-#1mu}\mkern#1mu}

\begin{document}

\begin{center}
\uppercase{\bf On Levine's Notorious Hat Puzzle}
\vskip 20pt

{\bf Joe Buhler}\\
{\smallit Reed College, Portland, OR, USA}\\
{\tt jpb@reed.edu}

\vskip 10pt

{\bf Chris Freiling}\\
{\smallit Project Inertia, San Diego, CA, USA}\\
{\tt chris.freiling@projectinertia.com}

\vskip 10pt

{\bf Ron Graham}\\
{\smallit UCSD, San Diego, CA, USA}

\vskip 10pt

{\bf Jonathan Kariv}\\
{\smallit Isazi Consulting, Johannesburg, South Africa}\\
{\tt jkariv@isaziconsulting.co.za}

\vskip 10pt

{\bf James R. Roche}\\
{\smallit Department of Defense, Fort Meade, MD, USA}\\
{\tt juggling.jim.roche@ieee.org}

\vskip 10pt

{\bf Mark Tiefenbruck}\\
{\smallit IDA, Center for Communications Research, La Jolla, CA, USA}\\
{\tt mgtiefe@ccrwest.org}

\vskip 10pt

{\bf Clint Van Alten}\\
{\smallit School of Computer Science and Applied Mathematics, University of the Witwatersrand, Johannesburg, South Africa}\\
{\tt clint.vanalten@wits.ac.za}

\vskip 10pt

{\bf Dmytro Yeroshkin}\\
{\smallit Géométrie Différentielle, Université Libre de Bruxelles, Brussels, Belgium}\\
{\tt Dmytro.Yeroshkin@ulb.ac.be}

\end{center}

\vskip 20pt
\centerline{\smallit Received: , Revised: , Accepted: , Published: }
\vskip 30pt

\pagestyle{myheadings}
\markright{\smalltt INTEGERS: 21 (2021)\hfill}
\thispagestyle{empty}
\baselineskip=12.875pt

\centerline{\bf Abstract}

\noindent
The Levine hat game requires $n$ players, each wearing an infinite random stack of black and white hats, to guess
the location of a black hat on their own head seeing only the hats worn by all the other players.
They are allowed a strategy session before the game, but no further communication. The players
collectively win if and only if all their guesses are correct.

In this paper we give an overview of what is known about strategies for this game, including an extended discussion of the case with $n = 2$ players (and a conjecture for an optimal strategy in this case). We also prove that
$V_n$, the optimal value of the joint success
probability in the $n$-player game, is a strictly decreasing function of $n$.

\vskip 30pt

\section{Introduction} \label{sec:intro}

In her blog in 2011, Tanya Khovanova \cite{TK} described
a ``hat puzzle'' from Lionel Levine 
involving $n$ people, each wearing a stack of infinitely many black and white hats;
she also gave a problem of her own inspired by that puzzle. 
Although superficially recreational, the Levine puzzle became notorious because of the difficulty of giving definitive answers 
to any of the questions it raised.

\smallskip

\noindent{\bf Levine Hat Puzzle:} A team with $n$ players has an initial strategy session, after which a 
referee places a stack of $h$ hats on each player's head.
Each hat is either black (a.k.a.~1) or white (a.k.a.~0).  The players must name a position on their own stack, 
and they collectively win if and only if they all
name the position of a black hat on their own head. 
Players cannot see the hats on their own heads, but they can see all of the other players' hats.
No communication between the players is allowed after the strategy session. 
Each of the $n \cdot h$ hats placed by the referee is chosen by an independent flip
of a fair coin (probability $1/2$ of each color). Each player $i$ 
must
communicate a positive integer $x_i$ to the referee (without 
learning
any of the values $x_j$ communicated to the
referee by the other players).
The players win if and only if for all $i$, the hat in position $x_i$ in the stack on the $i^{\rm th}$~player's head is black.
What (joint) strategy
should the players use to maximize their chance
of winning?

\smallskip

This puzzle seems to have arisen out of Levine's work with Tobias Friedrich
\cite{FL} on fast simulations of certain growth models. It is 
sometimes stated with wardens/prisoners or sultans/wise men instead of referees/players.

What can the players possibly do at the strategy session? 
They have to agree on a collection $\{f_i\}$ of ``strategy'' functions $f_i$, 
one for each player, that map the possible stacks that player $i$ might see to positive integers.  Each such 
(joint) strategy has a probability of success (based on the coin flips that will determine hat colors). 
Let $V_n^{(h)}$ denote the maximum value, over all possible strategies, of the probability of 
success for any given $n$ and $h$.  It is easy to see that that $V_n^{(h)}$ is nondecreasing as a function of~$h$. Let
\[
V_n := \lim_{h \to \infty} V_n^{(h)} = \sup_{h} V_n^{(h)} .
\]

Readers who like to work on puzzles themselves before seeing hints or solutions
are recommended to set this paper down immediately and prove
the following three
statements (which are easy, moderately challenging, and difficult, respectively): \\
(1) $V_n \ge 1/2^n$, ~~~ (2) $V_n \ge 1/(n+1)$, ~~~ (3) $V_n \ge c/\log(n)$ for some $c > 0$.

\medskip

The notoriety of the puzzle arises from the difficulty of answering the most basic
questions. In particular, no value of $V_n$ is known exactly for any $n > 1$,
and the limiting behavior of $V_n$ for large $n$
--- which was perhaps of primary interest to Levine ---
is unknown. He made the following conjecture.

\begin{conjecture}[Levine] \label{Conjecture:Levine}
The optimal success probability in the $n$-player game is $o(1)$
as $n$ goes to infinity; i.e.,
\[ \lim_{n \to \infty} V_n = 0. \]
\end{conjecture}

If you thought we were going to answer this question or find the value of
$V_2$, say, then you would, alas, be mistaken.
The goal of this paper is to describe several results that are aimed at these two big questions, in the hope that this will spur people to answer them.

Our main results are 
as follows:
(1) a proof that $V_n$ is strictly decreasing, i.e.,
\[
V_{n+1} < V_n \; \mbox{ for all } \; n;
\]
(2) the inequalities  
\[
\frac{7}{20} = 0.35 \; \le V_2 \; \le 0.3616\ldots
\]
(we conjecture that equality holds on the left); (3) a technique of ``matrix hints'' that, at least in principle, can be used to give 
arbitrarily good
upper bounds on the $V_n$; (4) an analysis of what happens when a fair coin is replaced by a Bernoulli coin that yields heads with probability~$p$; and (5) various ancillary results and data that might help someone who wants to answer any of the various open questions!
 
We learned about a fascinating recent preprint by Friedgut, Kalai, and Kindler~\cite{FKK} on the same day that they learned of ours; their paper also proves that $V_n$ is strictly decreasing, and conjectures generalizations that situate the problem in an interesting combinatorial and graph-theoretic context.

By way of introducing some of the basic techniques that will arise later, we now focus on the 
case $n = 2$.
Each of the players $A$ and $B$ (whom 
the reader
may think of as Alice
and Betty if that makes the problem seem more compelling)
has a large stack of $h$ hats on her head.
Asking
$A$ to say something about the random stack
of hats on her head, using only the information in the completely independent stack on $B$'s head, seems 
a bit unfair;
indeed, it feels like a mysterious game show, perhaps run by mathematicians with a strange sense of humor.
If $A$ and $B$ choose
random strategy functions $f_A$ and $f_B$, then each has 
an independent probability
of success $1/2$ of naming the position of a black hat on her head, so their joint probability of success is~$1/4$.
Of course, the trick is that they 
should jointly choose their strategy functions before the game 
so as to make their choices correlate in a useful way.

The following warm-up theorem proves weaker versions of the inequalities stated in (2) above, giving 
elementary precursors of ideas that will be developed later.
The lower bound is straightforward and has no doubt been found by many people who have looked at the puzzle.
The upper bound is trickier and has been discovered by (at least) several different
people. It seems likely that Noga Alon was the first; several of us first learned of it from
a letter that Walter Stromquist wrote to Ron Graham.

Before proving these bounds, 
we 
make two remarks that will be used
in the proof and 
will be assumed at many points later in the paper.

\begin{remark}\label{prob_vs_determ}
The reader might wonder whether the players
could do better with a probabilistic strategy.  
The success probability for any probabilistic
strategy is a convex linear combination of those for deterministic strategies, so (assuming
that the source of randomness is uncorrelated with the hat placements) there is 
always a deterministic strategy that is at least as good as any probabilistic strategy, and it suffices to consider deterministic
strategies throughout.  (This is 
also true
if the players receive ``hints'' as in the proof below).
\end{remark}

\begin{remark}\label{Vn_is_limit}
It is important to remember that $V_n$ is defined as a 
limit of success probabilities 
for
finite stacks of~$h$ hats.
For instance, in the proof of the next theorem 
there is a stack of $h$ hats on each player's head,
and the case in which 
some player has hats of only one color can
be ignored. This case has probability 
at most $2n / 2^{h}$
for any given $n$, which vanishes
in the limit as $h$ goes to infinity.
So we can, variously, speak of finite or infinite stacks of hats,
but must always remember that the infinite case is defined as a
limit of finite cases.
If strategies are actually allowed to use functions defined on
infinite sets, the situation is entirely different.
Any reader who can really see an infinite stack of hats all at once (as well
as perform computations on infinite sets in finite time), and who believes in the Axiom of Choice
and is not squeamish about non-measurable sets,
is advised to read Appendix A.
In that appendix, we consider the 1-person~(!) version of the game and
describe a strategy for which the player is ``virtually guaranteed" to win.
\end{remark}

\begin{theorem}\label{th1}
The optimal probability of success, $V_2$, in the 2-player game is bounded as follows:
\[
\frac{1}{3} \le V_2 \le \frac{3}{8}.
\]
\end{theorem}

\begin{proof}
The lower bound is proved by using the following ``first-black'' strategy:
Each player finds the position of the first (lowest) black hat on her
partner's head and names that 
same position on her own head!
For example, suppose that the hat stacks
begin as in the diagram below (with black hats represented by 1s and white hats by 0s).

\[
\begin{array}{ccc}
\vdots & & \vdots \\
0 & & 1 \\
1 & \rightarrow & 1 \\
0 & \leftarrow & 1 \\
0 & & 0 \\
  & &   \\
A & & B \\
\end{array}
\]

\bigskip

Positions are numbered from 1 starting at the bottom, as is fitting
for stacks of hats.
Then $A$ will say 2, $B$ will say 3, and they will lose.
(Although $B$ happens to have a black hat on level 3, $A$ does {\em not}
have a black hat on level 2.)

What is the probability that they will win using the first-1 (i.e., first-black) strategy? 
Consider the lowest level where
at least one of the players
has a 1. 
(By Remark 2, we may assume that there is such a level.)
There are 3 possible hat pairs for $A$ and $B$ at the critical level: 01, 10, and 11.
These are equally likely, so the first-black strategy has a 
success probability (or, as we 
sometimes say, ``value'') of 1/3.  
(Wow! 
Even though neither player has any information about the color of any hat on her own 
head, the team can do significantly better 
with correlated strategies than by making random guesses.)

Now we turn to the upper bound.
To prove it we use the curious device of allowing the players to get a hint from the 
referee. (This is a simple example of the ``matrix hints'' technique 
to be described later
for finding upper bounds.) This extra information certainly cannot lower
the optimal success probability of $A$ and $B$, since they are free to ignore the information if they wish.

Suppose that the referee takes pity on $A$ and $B$ and, before their strategy session,
shows the players a bit string $s$ (each bit having been determined by a flip of a fair coin) and
says that just before the game he will flip a fair coin one more time and then, based on
that flip, put the hat sequence corresponding either to $s$ or to its bitwise complement $s'$ 
on $A$'s head.
For instance, if $s$ is the bit string $101110\ldots $, then $A$ knows that her 
hat sequence will be either $s = 101110\ldots$
or $s' = 010001\ldots$, each with probability~$1/2$.  (Of course, after the game starts,
$B$ will actually see the sequence $s$ or $s'$ on $A$'s head.)
The players are given no information about the hat sequence on $B$'s head.

Note that if $A$ and $B$ refrain from looking at $s$, this game is exactly equivalent
to the originally described game: The referee is now determining $A$'s hat sequence 
in a two-stage process, but in the end, for any height $h$, all $2^h$ 
of her possible sequences 
are equally likely.
To establish the desired upper bound on $V_2$, we will show that for 
each possible 
value of $s$ and every possible joint strategy, the players' success 
probability is
at most $3/8$.

Because of the referee's hint, both players know during the strategy session
that $B$ will see one of
only two possible sequences on $A$'s head.
Thus $B$'s strategy is completely determined by the
integers $x = f_B(s)$ and $y = f_B(s')$ that
she will announce according to whether
she sees $s$ or~$s'$
as $A$'s stack. 
(It will turn out to be better for $B$ to choose $x \neq y$,
but the possibility that $x = y$ must also be considered.)
On the other hand, $A$'s strategy could be any function $f_A$ of the
stack that she sees on~$B$'s head.

For the rest of the proof, we suppose that a particular (though arbitrary) hint $s$
has been given. Then there are 2 possible values of $A$'s stack and $2^h$ possible
values of $B$'s stack, and all $2^{h+1}$ joint possibilities are equally likely. 
We also suppose that the team has chosen particular (though arbitrary)
strategy functions $f_A$ and $f_B$.

For each of the $2^h$ possible hat sequences $b \in \{0, 1\}^{h}$ that $B$ 
could be given, player $A$ will choose some level $f_A(b) \in \{1, 2, \ldots, h\}$
on her own head. Whatever level she chooses, her hat color at that level is determined
by a fair coin flip independent of the coin flips used to determine $s$ and $b$.
Thus we have the following ``atomic'' probability result:
\begin{equation} \label{eq:atomic_prob}
    \Prob(\mbox{$B$ has stack $b$; $A$ chooses a black hat}) = 1 / 2^{h+1}
     \mbox{~~for every $b$}.
\end{equation}

If $x = y$, then $B$ points to the same level whether she sees $s$ or $s'$ as $A$'s stack. 
Thus half of all values of $b$ (those with component $b_x = 0$) cause $B$ to choose a white hat,
in which case the team fails. The other $2^{h-1}$ values of $b$ (those with $b_x = 1$)
cause $B$ to choose a black hat, in which case the team wins if and only if $A$ chooses
a black hat. It follows from Equation~(\ref{eq:atomic_prob}) that
\begin{align*}
  & \Prob(\mbox{$A$ and $B$ both choose black hats})  \\
  & = \sum_{b: b_x = 1} \Prob(\mbox{$B$ has stack $b$; $A$ chooses a black hat})  \\
  & = (2^{h-1}) \cdot (1 / 2^{h+1}) = 1/4 < 3/8.
\end{align*}

If $x \neq y$, then $2^{h-2}$ values of $b$ have $b_x = b_y = 0$, in which
case $B$ will point to a white hat and the team will fail. For the 
remaining $3 \cdot 2^{h-2}$ values of $b$, the team has a chance. Again 
using Equation~(\ref{eq:atomic_prob}), we have
\begin{align*}
& \Prob(\mbox{$A$ and $B$ both choose black hats}) \\
& = \sum_{b: (b_x, b_y) \neq (0,0)} \Prob(\mbox{$B$ has stack $b$; $A$ and $B$ both choose black hats})  \\
& \leq \sum_{b: (b_x, b_y) \neq (0,0)} \Prob(\mbox{$B$ has stack $b$; $A$ chooses a black hat})  \\
& = (3 \cdot 2^{h-2}) \cdot (1 / 2^{h+1}) = 3/8,
\end{align*}
finishing the proof of the theorem.
\end{proof}

We invite the reader to verify the following two claims: If the levels
$x$ and $y$ are chosen in the above proof so that the hats in the hint $s$ at
those positions are of opposite colors, then
$A$ may be assumed without loss of generality to choose from those same 
two levels on her own head. Furthermore,
if the two players follow a ``first-black'' strategy within those two levels
when given the hint $s$,
they can achieve the upper bound 3/8 on their success 
probability in the limit as $h$ goes to infinity.

The main sections of this paper (a) consider strategies for 2 players
in detail (giving the lower bound $7/20$, and considering what happens
when the hat colors are determined by biased coins), (b) use ``matrix
hints'' to give upper bounds, including the $0.36 \ldots$ stated above, 
(c) give results for $n$ players, and (d) give lower bounds found
by computer for smallish~$n$.

Some of the results in this paper are recent, but many of the results and techniques here arose in
an extensive series
of emails in 2013-2015 between a somewhat amorphous group of people that included the authors as well as others; naturally, they referred to themselves as ``the Mad Hatters.''   
The origins of this paper lie in the Hatters'  desire
to collect the useful information in those much earlier emails.  
Many results, though not essential to the main flow of this paper, may be of 
interest or use to someone wanting to look more deeply into details, so they are included here as appendices.
Because preliminary drafts of the current paper were developed independently by at least two different subsets of authors, 
there are varying conventions (e.g., in the players' appellations and genders). We have, however, tried to 
maintain local consistency.

At the beginning of the process of revising, polishing, and extending the paper for the sake of publication, Ron Graham (1935-2020) left us. The Levine puzzle is the kind of question that delighted him, and he was fond of challenging people with  this 
particular
problem. In addition to pushing for greater clarity and more cleverness, he also repeatedly asked for more data. 
Partly at his urging, a rather large number of computational experiments were done on strategies for this puzzle in 2014, by many people.  
The other authors dedicate
this article to the memory of 
Ron's exuberance, mathematical and otherwise.

\section{New Results for Two Players}\label{Sec:2p}

\subsection{Proof that $V_2 \ge 7/20$}

We begin our discussion of two-player strategies with the promised $7/20$ lower bound.

\begin{theorem}\label{th7_20}
The optimal success probability $V_2$ satisfies
the inequality
\[
V_2 \; \ge \; \frac{7}{20} .
\]
\end{theorem}

\begin{proof}
To prove the lower bound, we exhibit a strategy with value $v = 0.35 = 7/20$.

Consider the following characterization of the first-black strategy. Each player looks at the first hat on the other player's head. If he does not like what he sees, then he skips it and looks at the next hat. The two players might not skip the same number of hats, but if they do, then they have an unusually high chance of winning. We can generalize this idea by letting each player look at the first $k$ hats before deciding whether to skip them. The simplest $7/20$ strategies use $k=3$.

Players $A$ and $B$
each look at the lowest three hats on their partner's head. They advance over a triple if it is all-0 or all-1 (or, as we will say, is ``monochromatic'') until they arrive at a non-monochromatic triple. Since there are 2 possible monochromatic triples and 6 non-monochromatic triples, the chance that a player will skip over a triple is $1/4$. They each will stop on a non-monochromatic triple, and then use the following algorithm: If there is a single 1 in the triple that they see, then they announce that the corresponding bit in their own string is also a~1. (So if 
$A$ skips twice and then the first bit in the next triple is a 1, 
$A$ would say~7.) If there are two 1-bits and one 0-bit in the lowest non-monochromatic triple that they see, then 
$A$
names the position {\bf A}bove the 0-bit, and 
$B$
names the position {\bf B}elow the 0-bit. Here above/below in the triple are interpreted cyclically; e.g., if 
$B$ skips one triple and then sees 011 (where 0, in the fourth position, is the lowest hat in that triple),
$B$ will say~6.

What is the value $v$ of this strategy? If $r$ denotes the probability that they win when both of their bottom triples are non-monochromatic, then
\[
v = \frac{1}{4} \cdot \frac{1}{4} \cdot v + \frac{1}{4} \cdot \frac{3}{4} \cdot \frac{1}{4} + \frac{3}{4} \cdot \frac{1}{4} \cdot \frac{1}{4} + \frac{3}{4} \cdot \frac{3}{4} \cdot r .
\]
The first term represents the case where they both see initial monochromatic triples (no harm, no foul, they just skip those and are then playing the same game). The second and third terms are the cases where one of 
$A$ or $B$ skips an initial triple but the other does not; this is the perhaps unfortunate case in which there can be no correlation because they are looking at different triples, so they are both making random uncorrelated guesses and have probability $1/4$ of winning. The fourth term represents the case where neither skips the initial triple. In order to solve this equation for~$v$, we have to calculate~$r$.

Call a non-monochromatic triple a {\it single} if it has a single 1 bit, and a {\it double} if it has exactly two 1 bits. There are 36 possible pairs of non-monochromatic triples: 9 cases where both are singles, 18 in which one is a single and one is a double (in one order or the other), and 9 in which both are doubles. Immediately below we show example pairs of all three types. (As it happens, all three examples are losing pairs for the above strategy.)
\[
\begin{array}{ccccccccccc}
0 && 0 & \hspace*{0.3in} 0 && 0 & \hspace*{0.3in} & 1 && 0 \\
0 &\leftrightarrow & 1 & \hspace*{0.3in} 1 & \leftrightarrow & 1 & \hspace*{0.3in} & 1 & \leftrightarrow & 1 \\
1 && 0 & \hspace*{0.3in} 0 && 1 & \hspace*{0.3in} & 0 && 1\\
 \\
A && B & \hspace*{0.3in} A && B & \hspace*{0.3in} & A && B
\end{array}
\]

In the 9 cases where they both have singles, they win exactly in the 3 cases where the 
1s are in the same location. In the 18 cases where they have a single and a double (for one or the other order), they are both correct only in the 6 cases where the 1 bit in the single is located in exactly the right location with respect to the two 1 bits in the double. In the case when 
$A$ and $B$ both have doubles, they fail {\it only} in the three cases where, as in the case pictured above, the bits are exactly aligned so that they are both wrong. This means that they are both correct in all 6 of the other cases, as the reader can check by trying the other two possibilities for 
$B$ when $A$
has 110 as pictured above. Thus 
\[
r = \frac{3+6+6}{36} = \frac{5}{12}.
\]
Solving the earlier equation for~$v$ gives $v = 7/20$ as claimed.
\end{proof}

The ``reset'' on black (not just on white) in the 7/20 strategy for the 2-person game was counterintuitive and surprisingly difficult to find; apparently it eluded discovery for three years after the puzzle was popularized in 2011 on Tanya Khovanova's blog \cite{TK}. 
It was finally found in 2014 by a California group (Larry Carter, Jay-C Reyes, Joel Rosenberg, 
and M.\ Tiefenbruck) and by a Pennsylvania group (J.\ Kariv and D.\ Yeroshkin).
The fact that the description above might be judged by some as ``easy to remember'' is probably a red herring, since a strategy amounts to nothing more than an arbitrary function from what a player sees to what they are supposed to do, and there is no reason that this needs to be structured in any memorable way whatsoever. In particular, there are multiple 7/20 strategies based on 3 hats (skipping monochromatic triples), and some are symmetric in the sense that both players have the same function.

It might seem paradoxical to skip over the case of 3 black hats, where one of the players is 
guaranteed to be right, but it seems to be necessary. One possible rationale is that it is 
harder to correlate with monochromatic triples, so it might be better just to skip them.
A more detailed explanation is given immediately below; it emphasizes why $A$ should reset when seeing 
a monochromatic triple, but the same argument applies when the roles of $A$ and $B$ are reversed.

If $B$ has a white triple, the team can win only if $B$ resets, so $A$ should hope for the best and reset as well.
(This argument is easy to find.) If $B$ has a black triple and $A$ has a non-monochromatic triple, then $B$ is not
going to reset in any case, so the team will will on average win half the time whether or not $A$ resets.
Resetting upon seeing a black triple 
thus makes a net difference only when {\em both} players have monochromatic triples, at least 
one of them black. (As we will see when considering the $n$-player game, the actual requirement 
is that at {\em most} one of these monochromatic ``tiers" be {\em white}.) Within this subset 
of 3 cases, a shared strategy to reset on seeing a black triple (as well as on seeing a white 
triple, which is assumed) gives up 1 sure win and 2 sure losses in order to get 3 fresh starts 
at a new 3-level tier. As long as the team has a basic (non-resetting) 3-level strategy that wins 
strictly more than 1/3 of the time (it is 22/64 for the 2-player game), incorporating the reset 
on seeing a black triple is a net win.

One might expect that strategy functions based on larger numbers of hats would yield
increased probability of success.
We do not know for sure, but 
we think otherwise and have
the vague intuition that with larger clumps it is harder to usefully correlate assignments of probability mass to various choices. Many attempts were made in 2014 to find better strategies, but none succeeded. 
On these grounds,
we have come to believe that the 7/20 strategies are quite possibly optimal.

\begin{conjecture}
We conjecture that
\[
V_2 \; = \; \frac{7}{20}.
\]
\end{conjecture}

\subsection{The Two-Player Game with a Biased Coin}

The $7/20$ strategy
for the 2-player Levine puzzle has been a sticking point for about seven years. No one has found a better strategy or proved that it is optimal. Sometimes it is useful to change a problem when stuck, and one natural idea here is to replace the flip of a fair coin, used by the referee to determine hat color, by a biased coin which has probability $p$ of giving a black hat. This turns out to have several virtues, one of which is to show that there are actually several distinct optimal strategies for $p = 1/2$ that are {\em not} equivalent to each other when $p \neq 1/2$. As is usual, we set $q=1-p$.
 
The special case with $p=\frac{a}{b}$ rational 
(for some $b > 2$) corresponds to the natural extension of having 
$b$ hat colors, of which some ($a$) are considered good  and the remaining ones are considered bad. The goal is then for each player to choose a hat of a good color (equivalently, for none of them to choose a bad color). The special cases of $a=1$ and $a=b-1$ are of particular interest as they respectively correspond to the case of a single good color and a single bad (nuclear) color.
 
In this section we will consider general values of $p \in (0,1)$ and compare the team's success probability, or value, using different strategies $S$. Thus we extend the notation from the previous section to consider quantities
\[
V_n^{(h)}(p; S).
\]
We will omit the superscript $h$ when considering the limiting case as $h \rightarrow \infty$ and will omit the argument $S$ when referring to the supremum over all strategies. If the argument $p$ is also omitted, then the default value $p = 1/2$ is understood.

Using the naive strategy of each player choosing the first hat corresponding to a black hat on 
the other player's head, we can obtain a probability of winning of 
$\frac{p}{2-p}$. If each player instead chooses the first hat corresponding to a {\em white} hat on the partner's head, the probability of winning is
$\frac{2p^2}{1+p}$.

As for the special case of $p=\frac{1}{2}$, these strategies, while easy to state and better than random, are not optimal. We construct four distinct strategies based upon the 3-hat strategy that all achieve a performance of $0.35$ for the special case of $p=0.5$. However they all perform differently for general $p$. 

\begin{theorem}\label{Thm:Distinct}
There exist at least four distinct strategies $S$ that achieve $V_2(1/2; S) = 7/20$ but are inequivalent for $V_2(p; S)$ when $p_2 \neq 1/2$.
\end{theorem}

We shall construct the four strategies mentioned above and will denote them by $S_1$, $S_2$, $S_3$ and $S_4$. We provide their respective win rates here 
for the reader's convenience:
\begin{align*}
V_2(p; S_1) &= \frac{p (1 + p + p^2 + 3 p^3 - 3 p^4 + p^5)}{2 + p + p^2 + p^3 - p^4} ;\\
V_2(p; S_2) &= \frac{p(1-p+p^2+p^3)}{2-3p+3p^2} ;\\
V_2(p; S_3) &= \frac{p(1 + 5p - 10p^2 + 10p^3 - 5p^4 + p^5)}{(2-2p+p^2)(1+p)(2-p)} ;\\
V_2(p; S_4) &= \tfrac{p ( 1 + 7p - 21p^2 + 35p^3 - 20p^4 - 14p^5 + 40p^6 - 48p^7 + 40p^8 - 22p^9 + 7p^{10} - p^{11} )}{(1-p+p^2)(1+p-p^2)(2-2p+p^2)(1+p^2)(1+p)(2-p)}.
\end{align*}

These combine to give a lower bound for $V_2(p)$ of 
\begin{enumerate}
\item $\displaystyle \frac{p (1 + p + p^2 + 3 p^3 - 3 p^4 + p^5)}{(1+p)(2-p)(1+p^2)} \leq V_2(p)$ for $p \leq \frac{1}{2}$;

\item $\displaystyle \frac{p(1 + 5p - 10p^2 + 10p^3 - 5p^4 + p^5)}{(2-2p+p^2)(1+p)(2-p)} \leq V_2(p)$ for $\frac{1}{2} \leq p$.
\end{enumerate}

The curve given by this theorem is provided in Figure \ref{Fig:Bounds}, in the next section.

\subsection{Constructing Strategies}\label{Sec:Strats}

The strategies described in this 
subsection are based on ones found by a computer search for the game with only finitely many hats on each player's head. In particular, an exhaustive search was run to find the optimal strategy with three hats and the optimal ``symmetric'' strategy (i.e., both players use the same strategy) with four hats. Beyond these two cases, the authors ran hill-climbing and genetic algorithms for up to 12 hats; no strategies were found with better performance than the ones described here. In Section~\ref{sec:3Hats} we describe the outcome of the search of strategies with only three hats, and then in Section~\ref{sec:InfHats} we adapt the results to obtain the best known strategies for infinitely many hats. For related work on applying genetic algorithms to hat problems, see \cite{BGK}, which considers a different hat game.

\subsubsection{Basic Strategy for Three Hat Levels}\label{sec:3Hats}

\begin{table}[hb]
\centering
\begin{tabular}{|c||c|c|c|c|c|c|c|c|}
\hline
Black hats & $\emptyset$ & $\{1\}$ & $\{2\}$ & $\{1,2\}$ & $\{3\}$ & $\{1,3\}$ & $\{2,3\}$ & $\{1,2,3\}$\\
\hline
Picture & $\square\square\square$ & $\blacksquare\square\square$ & $\square\blacksquare\square$ & $\blacksquare\blacksquare\square$ & $\square\square\blacksquare$ & $\blacksquare\square\blacksquare$ & $\square\blacksquare\blacksquare$ & $\blacksquare\blacksquare\blacksquare$\\
\hline
Choice & any & 1 & 3 & 1 & 2 & 2 & 3 & any\\
\hline
\end{tabular}
\vspace{0.2cm}
\caption{$S_0$, an optimal strategy on 3 hats}\label{Table:3Hats}
\end{table}

We denote by $S_0$ the symmetric strategy 
for 3 hat levels defined by Table~\ref{Table:3Hats}. 
As shown in Table~\ref{Table:3HatsWin},
strategy $S_0$ wins in $22$ of $64$ cases. The likelihood of each case depends on $p$ as is also shown in Table~\ref{Table:3HatsWin}. For convenience, we assume that a player who sees all hats of the same color points to the first hat on his or her own head. The columns correspond to the distribution of black hats on 
the first player's head
and the rows to the distribution on 
the second player's.
The cells are blank when the players lose. It is useful to note that this strategy is the unique optimal strategy for every value of $p$, up to reordering the hats on one or both of the players' heads.

\begin{table}[hb]
\centering
\begin{tabular}{|c||c|c|c|c|c|c|c|c|}
\hline
& $\square\square\square$ & $\blacksquare\square\square$ & $\square\blacksquare\square$ & $\blacksquare\blacksquare\square$ & $\square\square\blacksquare$ & $\blacksquare\square\blacksquare$ & $\square\blacksquare\blacksquare$ & $\blacksquare\blacksquare\blacksquare$\\
\hline\hline
$\square\square\square$ & & & & & & & &\\
\hline
$\blacksquare\square\square$ & & $p^2q^4$ & & $p^3q^3$ & & & & $p^4q^2$ \\
\hline
$\square\blacksquare\square$ & & & & & $p^2q^4$ & $p^3q^3$ & & \\
\hline
$\blacksquare\blacksquare\square$ & & $p^3q^3$ & & $p^4q^2$ & & $p^4q^2$ & & $p^5q$ \\
\hline
$\square\square\blacksquare$ & & & $p^2q^4$ & & & & $p^3q^3$ &\\
\hline
$\blacksquare\square\blacksquare$ & & & $p^3q^3$ & $p^4q^2$ & & & $p^4q^2$ & $p^5q$ \\
\hline
$\square\blacksquare\blacksquare$ & & & & & $p^3q^3$ & $p^4q^2$ & $p^4q^2$ & \\
\hline
$\blacksquare\blacksquare\blacksquare$ & & $p^4q^2$ & & $p^5q$ & & $p^5q$ & & $p^6$ \\
\hline
\end{tabular}

\caption{Winning combinations for $S_0$ with probability of each event given in the block}\label{Table:3HatsWin}
\end{table}

The sum of all the winning probabilities is as follows:
\[
V_2^{(3)}(p; S_0) = 3p^2q^4 + 6p^3q^3 + 8p^4q^2 + 4p^5q + p^6 = 3p^2 - 6p^3 + 8p^4 - 6p^5 + 2p^6.
\]

\subsubsection{Adaptation to Infinitely Many Hat Levels}\label{sec:InfHats}

We give four adaptations of the above 3-hat strategy to the general game that performed well in cases of up to 12 hats for various values for $p$. 
As explained below, we assume without loss of generality that the first 3 hats considered are
those in positions 1 to 3, that the second group of 3 hats occupies either 
positions 3 to 5 or 4 to 6, and so on.

\smallskip

Strategy $S_1$:
\begin{enumerate}
\item If the first three hats of the other player are not monochromatic, play the 3-hat strategy $S_0$.

\item If the first three hats of the other player are BBB or WWW, repeat $S_1$ on hats 3 through $\infty$.
\end{enumerate}

Strategy $S_2$:
\begin{enumerate}
\item If the first three hats of the other player are not monochromatic, play the 3-hat strategy $S_0$.

\item If the first three hats of the other player are BBB or WWW, repeat $S_2$ on hats 4 through $\infty$.
\end{enumerate}

The strategies $S_3$ and $S_4$ are constructed the same way as $S_1$ with different, but equivalent, symmetric 3-hat strategies. We provide those 3-hat strategies in Tables~\ref{Table:S3}~and~\ref{Table:S4}.

\begin{table}[hb]
\centering
\begin{tabular}{|c||c|c|c|c|c|c|c|c|}
\hline
Black hats & $\{1\}$ & $\{2\}$ & $\{1,2\}$ & $\{3\}$ & $\{1,3\}$ & $\{2,3\}$\\
\hline
Picture &$\blacksquare\square\square$ & $\square\blacksquare\square$ & $\blacksquare\blacksquare\square$ & $\square\square\blacksquare$ & $\blacksquare\square\blacksquare$ & $\square\blacksquare\blacksquare$\\
\hline
Choice & 3 & 2 & 2 & 1 & 3 & 1\\
\hline
\end{tabular}
\vspace{0.2cm}
\caption{3-hat strategy that produces $S_3$}\label{Table:S3}
\end{table}

\begin{table}[hb]
\centering
\begin{tabular}{|c||c|c|c|c|c|c|}
\hline
Black hats & $\{1\}$ & $\{2\}$ & $\{1,2\}$ & $\{3\}$ & $\{1,3\}$ & $\{2,3\}$\\
\hline
Picture & $\blacksquare\square\square$ & $\square\blacksquare\square$ & $\blacksquare\blacksquare\square$ & $\square\square\blacksquare$ & $\blacksquare\square\blacksquare$ & $\square\blacksquare\blacksquare$\\
\hline
Choice & 2 & 1 & 1 & 3 & 2 & 3\\
\hline
\end{tabular}
\vspace{0.2cm}
\caption{3-hat strategy that produces $S_4$}\label{Table:S4}
\end{table}

\begin{remark}
Iterating over all the basic 3-hat strategies equivalent to $S_0$, if one uses the shift-by-3-levels construction for $S_2$, the strategies will once again be equivalent (this includes the asymmetric strategy described in Theorem~\ref{th7_20}). On the other hand, if one undertakes the shift-by-2-levels construction for $S_1$, then the win rate of the strategy will be one of $V_2(p;S_1)$, $V_2(p;S_3)$, or $V_2(p;S_4)$. 
\end{remark}

We close this section by computing the success probability (or value) for the $S_2$ strategy, which has the most concise presentation. The computation for $S_1$, which yields our best lower bound for $p \leq 1/2$, is more complicated and is deferred to Appendix~\ref{Sec:S1App}. The computations for $S_3$ and $S_4$ are omitted, since they follow the same procedure as that for $S_1$.

\subsection{Computing the Performance of $S_2$}

Consider the following cases:

\begin{enumerate}
\item[(a)] Neither player has WWW or BBB as their first three hats. The probability of this occurring and 
the players' winning
is $3p^2q^4 + 6p^3q^3 + 6p^4q^2$ (see Table~\ref{Table:3HatsWin}).

\item[(b)] One player has BBB as their first three hats, and the other does not have either WWW or BBB. 
Either of the two players can have the BBB stack, which occurs with probability of $p^3$, and the probability of not having WWW or BBB is $1 - p^3 - q^3$. The probability of this case occurring is therefore $2p^3(1 - p^3 - q^3)$. The probability of winning given this case is $p$ since the player with BBB guesses correctly, and the other player chooses a hat in position 4 or greater, which has probability $p$ of being B. Thus, $2p^3(1 - p^3 -q^3)p$ is the probability of this case occurring and 
the players' winning.

\item[(c)] One player has WWW as their first three hats, and the other does not have either WWW or BBB. In this case the players are certain to lose.

\item[(d)] Both players have either BBB or WWW as their first three hats, which occurs with probability $p^6 + 2p^3q^3 + q^6 = (p^3 + q^3)^2$. Then the strategy looks to the next three hats, and this repeats, giving an infinite sum with ratio $(p^3 + q^3)^2$.
\end{enumerate}

The overall probability of a win is therefore
\begin{align*}
V_2(p; S_2) & = \frac{(3p^2q^4 + 6p^3q^3 + 6p^4q^2) + 2p^3(1 - p^3 -q^3)p}{1 - (p^3 + q^3)^2} \\
& = \frac{p(1-p+p^2+p^3)}{2-3p+3p^2}.
\end{align*}

\section{The Matrix Game and Upper Bounds for Two Players}\label{Sec:Matrix}

\subsection{An Informal Introduction to the Matrix Game}

In the introductory section of this paper, we considered giving a simple hint to the two players. This allowed us to compute an upper bound of $3/8$ for the value of the two-person game where the probability of each hat color is $1/2$. In 
the current section we present a scheme to construct more elaborate hints that prove better upper bounds.
For variety, we refer to the players as Alice and Bob in the 
present subsection.
In the later subsections of Section~\ref{Sec:Matrix}, we will use the
more rigorous-sounding names ``Player 1'' and ``Player 2'' but will 
continue to think of the two players as female and male, respectively.

Recall that in Theorem~\ref{th1} in  
the introduction,
the referee helped the players by revealing a bit string that was guaranteed to be either the first player's exact 
hat sequence or its bitwise complement.
One way for the referee to produce a pair of complementary sequences is to randomly choose rows from the $2\times 2$ identity matrix,
\[
\begin{pmatrix}
1 & 0 \\
0 & 1
\end{pmatrix}.
\] 
Note that in order to conform to the placing of hats on heads, all rows in this section will be indexed from the bottom up, starting with index 1. For example, choosing row 2, row 1, row 1, row 2, row 1, ... would produce the following complementary pair of hat sequences:
\[
\begin{pmatrix}
\vdots & \vdots\\
0 & 1 \\
1 & 0 \\
0 & 1 \\
0 & 1 \\
1 & 0
\end{pmatrix}.
\]
To construct a more elaborate hint, we could repeat this procedure using a different matrix, perhaps
\[
\begin{pmatrix}
0 & 0 & 1 & 0 & 1 & 1 \\
0 & 1 & 0 & 1 & 0 & 1 \\
1 & 0 & 0 & 1 & 1 & 0 \\
1 & 1 & 1 & 0 & 0 & 0 
\end{pmatrix}.
\]
We hope that a larger matrix will amount to a weaker hint given by the referee, which might give us a tighter upper bound. Choosing rows randomly gives us a collection of six 
infinite hat sequences.
For example,
\[
\begin{pmatrix}
\vdots \\
x_4 \\
x_3 \\
x_2 \\
x_1 
\end{pmatrix}
\in
\begin{pmatrix}
\vdots & \vdots & \vdots & \vdots & \vdots & \vdots & \\
0 & 0 & 1 & 0 & 1 & 1 \\
1 & 0 & 0 & 1 & 1 & 0 \\
0 & 1 & 0 & 1 & 0 & 1 \\
1 & 0 & 0 & 1 & 1 & 0 
\end{pmatrix}
\]
would be produced if the referee chose row 2, row 3, row 2, row 4, .... The hint comes when the referee reveals these six sequences and promises that 
Alice's sequence 
is among them. 

Notice that in the example above, hats $x_1$ and $x_3$ are guaranteed to be identical, no matter which of the six sequences is chosen. More generally, when 
Alice chooses an index, it only matters which row was used to produce the hats at that level, and this row is known to all. Thus, once the matrix has been fixed, the game on 
Alice's head is finite.
The referee's task is merely to choose a random column of the matrix,
and Alice's decision is reduced to identifying a row of the matrix. For the players to have a chance of winning, there must be a 1 in the resulting row and column.

In the case of this $4\times 6$ hint, it may be convenient to imagine that Alice has a single hat in the shape of a tetrahedron. The referee chooses an edge on the hat, and Alice hopes to choose a vertex adjacent to that edge.

So, the matrix hint greatly simplifies the game for Alice. But what about Bob, who receives no such hint? Bob will see which of the six columns is randomly chosen by the referee and will base his decision on this observation. Assuming a deterministic strategy, there are at most six hats that Bob will ever use. It may be fewer than six, because Bob may decide to use the same hat index for several different columns. In fact, Bob's strategy boils down to choosing some partition of the six columns of the matrix. So, for example, if the six columns are $c_1,c_2,c_3,c_4,c_5,c_6$, then Bob may decide on the partition $\{ \{c_1,c_4,c_5\}, \{c_2,c_3\}, \{c_6\} \}$, which means that he will choose a certain hat --- which may as well be $x_1$ --- when he sees $c_1$, $c_4$, or $c_5$; choose $x_2$ when he sees $c_2$ or $c_3$; and choose $x_3$ when he sees $c_6$. So the game for Bob also reduces to choosing from among a fixed, finite set of strategies.

Once Bob's strategy is fixed, there is an obvious best strategy for Alice, which we now describe. Alice observes the hats on Bob's head. Knowing Bob's strategy, Alice can see which of the six columns would cause Bob to succeed. Call these ``winning columns." In order to survive, they need the referee to choose one of these winning columns, resulting in a 1 on Bob's head. But they also need a 1 on Alice's head, so they also need the referee to choose a column with a 1 in Alice's chosen row, whatever that may be. The probability of winning, then, is $k/6$, where $k$ is the number of winning columns with a 1 in Alice's row. So Alice simply chooses any row that maximizes this probability. Finding the expected value of this probability over all of Bob's hat assignments gives the value of Bob's strategy.

For a concrete example, suppose once more that Bob chooses the partition $\{ \{c_1,c_4,c_5\}, \{c_2,c_3\}, \{c_6\} \}$. 
Then Bob is using a 3-hat strategy. We consider the eight possible assignments of these three hat colors, shown in Table~\ref{tab:table1}. So the value of this strategy is $17/48$.

\begin{table}[htb]
\centering
\begin{tabular}{l|c|c|r}
\textbf{Bob's colors} & \textbf{Winning Columns}& \textbf{Best row for Alice} & $\Prob(\mathrm{win})$\\
\hline
000 & $\{\}$ & does not matter & $0/6$ \\
001 & $\{c_6\}$ & row 3 & $1/6$ \\
010 & $\{c_2,c_3\}$ & row 1 & $2/6$ \\
011 & $\{c_2,c_3,c_6\}$ & row 1 & $2/6$ \\
100 & $\{c_1,c_4,c_5\}$ & row 2 & $3/6$ \\
101 & $\{c_1,c_4,c_5,c_6\}$ & row 2 & $3/6$ \\
110 & $\{c_1,c_2,c_3,c_4,c_5\}$ & row 1 & $3/6$ \\
111 & $\{c_1,c_2,c_3,c_4,c_5,c_6\}$ & does not matter & $3/6$
\end{tabular}
\caption{Computing Alice's strategy}
\label{tab:table1}
\end{table}

To find the best overall strategy, loop over all partitions of the columns of the matrix. For each of these ``Bob strategies," compute its value. Then choose the partition with the best value. This maximum value gives an upper bound for the two-player game.

One might think that after all this work, the $4 \times 6$ hint would give an improvement on our $3/8$ upper 
bound. In that case, one would be wrong! 
All we get is another way to prove the $3/8$ bound. However, applying an even more elaborate $8\times 14$ hint,
\setcounter{MaxMatrixCols}{20}
\[
\begin{pmatrix}
0 & 0 & 0 & 0 & 0 & 0 & 0 & 1 & 1 & 1 & 1 & 1 & 1 & 1\\
0 & 0 & 1 & 0 & 1 & 1 & 1 & 0 & 0 & 0 & 1 & 0 & 1 & 1\\
0 & 1 & 0 & 1 & 0 & 1 & 1 & 0 & 0 & 1 & 0 & 1 & 0 & 1\\
0 & 1 & 1 & 1 & 1 & 0 & 0 & 1 & 1 & 0 & 0 & 0 & 0 & 1\\
1 & 0 & 0 & 1 & 1 & 0 & 1 & 0 & 1 & 0 & 0 & 1 & 1 & 0\\
1 & 0 & 1 & 1 & 0 & 1 & 0 & 1 & 0 & 1 & 0 & 0 & 1 & 0\\
1 & 1 & 0 & 0 & 1 & 1 & 0 & 1 & 0 & 0 & 1 & 1 & 0 & 0\\
1 & 1 & 1 & 0 & 0 & 0 & 1 & 0 & 1 & 1 & 1 & 0 & 0 & 0 
\end{pmatrix},
\]
does finally give an improvement of $ 81/224 = 0.361607\ldots$. It is worth mentioning that there are 3920 different partitions that provide this bound, which are of 8 different types once we account for equivalences under permutations of rows and columns. 

There are a couple of things to keep in mind when constructing hints. First, the players are free to ignore the hint if they wish. Second, for a matrix hint each column in the matrix must have equal numbers of zeros and ones. Otherwise, the referee will not fulfill her obligation that the final sequence be random with $p=1/2$. Third, the stronger the hint, the weaker the upper bound will be. So, to produce good bounds, we would like a matrix that is short (few rows) and wide (many columns). Also, we would like to give such hints to as few players as possible. The catch is that we loop over all partitions of the matrix columns, and this quickly becomes infeasible as the number of columns is increased.

The $8\times 14$ hint was constructed by considering the fourteen non-constant affine functions on three bits. The $81/224$ upper bound it produces remains the best provable upper bound to date for $V_2$. A reader interested in beating this record may be tempted to try the shorter and wider $6\times 20$ matrix formed by the 20 three-element subsets of six elements. Surprisingly, this matrix does not do as well. Although an exhaustive search over column partitions was not performed, strategies achieving $117/320 = 0.365625$ have already been found, demonstrating that the upper bound will be worse. It is possible that a $16 \times 30$ matrix derived from non-constant affine functions on four bits might yield a better upper bound, since hill climbs have found no success rate higher than $0.35625$ for this larger matrix. However, an exhaustive search over all partitions of the 30 columns of this matrix currently seems infeasible.

The rest of this section is dedicated to formalizing these hints, using them to compute some bounds, applying them to other values of $p$, and providing a proof that these upper bounds converge to the optimal success probability.

\subsection{The Matrix Game}

We now turn our attention to computing upper bounds on $V_2(p)$. We do this by introducing a matrix game, which will be isomorphic to the hats game with the players given a little extra information. This extra information allows us to compute upper bounds. Furthermore, we shall show that by choosing a sufficiently large matrix, we can make these upper bounds arbitrarily tight. 

\medskip

\noindent
{\bf Note:} Throughout most of the paper, the variable $n$ is used for the number of players. In this section, however, the number of players is always 2, and we use $n$ to refer to the number of columns in a matrix.

\medskip

Given a Bernoulli parameter $p \in (0,1)$ and a matrix of $\{0,1\}$-valued entries,
\[
M = 
\begin{pmatrix}
x_{1,1} & x_{1,2} & \cdots & x_{1,n} \\
x_{2,1} & x_{2,2} & \cdots & x_{2,n} \\
\vdots & \vdots & \ddots & \vdots \\
x_{m,1} & x_{m,2} & \cdots & x_{m,n} 
\end{pmatrix},
\]
we can play the following two-player cooperative game. 

\textbf{Player 2:} 
The second player chooses 
an equivalence relation, $\sim$, on the set $\{1, \dots, n\}$ (or, equivalently, a partition of the set of columns of $M$). Each possible choice for $\sim$ yields $C$ equivalence classes (or disjoint subsets of the $n$ columns) for some $C$ with $1 \le C \le n$. (This equivalence relation corresponds to how Player 2 will choose a hat level on his own head based on which of $n$ ``hint columns'' the referee ultimately selects for the stack of hats on Player 1's head.)

\textbf{Referee:} 
The referee independently
flips a Bernoulli($p$) coin $C$ times, once for each of the $C$ column subsets chosen by Player 2, assigning the value 1 with probability $p$ and the value 0 with probability $1-p$ each time. (For the hat game, these correspond to the hat colors on the $C$ distinguished levels from which Player 2 will choose a hat on his own head.)

The referee then forms a vector $v = (v_1, \dots, v_n)$, where each entry $v_j$ is the binary value assigned above to the component of the partition that contains column $j$. We say that the vector $v$ {\em respects} the equivalence relation $\sim$, because equivalent columns are assigned the same binary value.

\textbf{Player 1:} 
The first player observes
the now-specified vector $v$ assigned by the referee and chooses a row, $r$, of $M$ and reports the dot product $r \cdot v$.

\medskip

Players 1 and 2 wish to maximize the dot product. As Player 1 can easily compute the dot product for each possible row $r$, it is trivial for her to choose the maximal dot product for any given $v$.

\medskip

\textbf{Conditional Value:} The conditional value $V(M; p, \sim)$ of the matrix game on $M$ --- abbreviated to $V(M; \sim)$ when $p$ is understood --- for a particular equivalence relation is the expected value (over all possible realizations of the referee's Bernoulli($p$) coin flips) of the reported (maximal) dot product $r \cdot v$ divided by $n$, the number of columns. That is,
\[
V(M; p, \sim) = \frac{1}{n} \underset{v} \sum \underset i \max (r_i \cdot v) \Prob(v), \\
\] 
where the probability of each vector, $\Prob(v)$, depends on the value of $p$ and on the particular choice of $\sim$. (The normalization by $1/n$ corresponds in the hats-with-hints game to the fact that the referee will create $n$ possible stacks for Player 1's head and reveal this set of $n$ possibilities to both players the night before the game; on the day of the game, the referee will uniformly at random choose one of these $n$ possible stacks, which Player 2 will see but Player 1 will not.)

\medskip

\textbf{Value:} The (unconditional) value of the matrix game $V(M; p)$ --- abbreviated to $V(M)$ when $p$ is understood ---
is the maximum of the conditional value of the matrix game over all choices of equivalence relation:
\[
V(M; p) = \underset{\sim} \max \,V(M; p, \sim ).
\]

The matrix game was devised primarily as a way of formalizing the hints that the referee can give the players to yield an upper bound on the value of the original hat game. However, by using matrices with slightly different characteristics, we can also encode the original hat game {\em without} hints and derive {\em lower} bounds on the value of the original hat game.

Theorem~\ref{Thm:UpperLower} asserts that these two bounds converge asymptotically, so that we could, in principle, approximate $V_2(p)$ arbitrarily closely (at least for rational $p$) by choosing suitable large matrices.

\subsection{Lower Bounds from the Matrix Game}

For a fixed rational value of $p$, we can get a lower bound on the probability of winning the original hat game by playing the matrix game on a matrix with appropriately repeated columns, such as
\[
L_{3,\frac{1}{2}} = 
\begin{pmatrix}
0 & 0 & 0 & 0 & 1 & 1 & 1 & 1 \\
0 & 0 & 0 & 0 & 1 & 1 & 1 & 1 \\
0 & 1 & 0 & 1 & 0 & 1 & 0 & 1 
\end{pmatrix}
\]
or 
\[
L_{2,\frac{2}{3}} = 
\begin{pmatrix}
0 & 0 & 0 & 1 & 1 & 1 & 1 & 1 & 1 \\
0 & 1 & 1 & 0 & 0 & 1 & 1 & 1 & 1 \\
\end{pmatrix}.
\]
This is equivalent to playing the hat game where Player 2 looks only at the first few hats of Player 1, but Player 1 can look at all of Player 2's hats. The columns of $L_{3,\frac{1}{2}}$ represent the possible colorings of Player 1's first three hats with $p=\frac{1}{2}$, while $L_{2,\frac{2}{3}}$ represents the possible colorings of her first 2 hats with $p=\frac{2}{3}$ with the repeated columns representing the proportionate likelihood of the colorings. Player 2's hat strategy assigns one of his own hat positions to each distinct coloring. This gives an equivalence relation on the columns, where two columns are equivalent if they are assigned to the same hat position on Player 2's head. (Repeated instances of the same column may without loss of generality be assigned to the same equivalence class, corresponding to the fact that the players' strategies for the hat game may be taken to be deterministic.)

A vector $v$ corresponds to a coloring on Player 2's head; $v_i$ is the value of the hat in the position on Player 2's head chosen by the given strategy. When Player 1 chooses the best row, she is really choosing a hat on her own head that is most likely to give a pair of matching ones when Player 2 uses his strategy.

The conditional value $V(M; \sim)$ is the probability of both players guessing black hats for a given choice of equivalence relation, and the value $V(M)$ is the probability of this occurring for an optimal strategy.

This gives a lower bound for the complete hat game because it restricts the hats that Player 2 can look at, but does not give any advantage over the original game.

\subsection{Upper Bounds from the Matrix Game}

To get an upper bound with $p = \frac{a}{b}$, we take any matrix of 
0s and 1s in which the proportion of 
1s is $p$ in every column. The matrix can have duplicate columns. Then, when we play the game on this matrix, we get an upper bound. As an illustration, consider the following matrix:
\[
U = 
\begin{pmatrix}
0 & 0 & 1 & 1 & 1 \\
1 & 1 & 1 & 0 & 1 \\
1 & 1 & 0 & 1 & 0 
\end{pmatrix}.
\]
We randomly generate an infinite sequence of rows from $U$. This produces,
in this case, five infinite (vertical) sequences of bits.
Since exactly two thirds of the bits in each column are ones,
each of the five infinite sequences has each bit independently equal to 1 with probability $p = \frac{2}{3}$. We then reveal to the two players that Player 1's hat sequence will be randomly chosen from this set of five. Although the five infinite sequences are not independent, this final selection will still be a randomly chosen sequence of hats. Player 2's hat strategy will assign one of his own hat positions to each of these 5 possible sequences for Player 1. As before, this creates an equivalence relation on the five columns of $U$, where two columns are equivalent if the sequences they generated are assigned by Player 2 to the same hat position on his own head. Since we generated infinitely many independent selections from the rows of $U$, by the second Borel-Cantelli lemma each row is almost surely chosen infinitely often. Therefore, when Player 1 chooses one of her own hats, it is equivalent to choosing one of the rows of $U$. The reason that this is an upper bound for the hat game is because 
the players are given extra information. They are not restricted in any way and they do not have to use this extra information if they do not want to, so this cannot hurt them. But it may help.

\subsection{The Meeting of Upper and Lower Bounds}

For each rational probability $p=a/b$ in lowest terms and each positive integer $m$ that is a multiple of $b$, we define two matrices, $L_{m,p}$ and $U_{m,p}$. The columns of these matrices will be elements of $\{0,1\}^m$. In $L_{m,p}$ all $2^m$ such columns appear, and each column with $t$ 
1s occurs $a^t (b-a)^{m-t}$ times, for a total of $b^m$ columns. In $U_{m,p}$ only the $\binom{m}{mp}$ columns with $mp$ 
1s occur, and there is no repetition of columns. For both matrices, the columns may be ordered arbitrarily.

As argued in the previous sections, $V(L_{m,p})$ is a lower bound for the value of the two-person hat game with black-hat probability $p$, and $V(U_{m,p})$ is an upper bound. We now state the main theorem of this section, that the matrix-based upper and lower bounds converge to $V_2(p)$, the value of the 2-player hat game with (rational) black-hat probability $p \in (0,1)$.

\begin{theorem}[Convergence Theorem] \label{Thm:UpperLower}
Let $L_{m,p}$ and $U_{m,p}$ be defined as above. Then
\[
\lim_{m \rightarrow \infty} V(U_{m,p}) - V(L_{m,p}) = 0.
\]
\end{theorem}

We prove this result in Appendix~\ref{Sec:MatrixBounds}.

For certain values of $p$ (in particular, those of the form $\frac{1}{k}$ or $\frac{k-1}{k}$), the general techniques above can be applied to particular matrices that yield rather good upper bounds without requiring too much work. In Appendix~\ref{App:SpecificUpperBounds} we prove the following general theorem and display hint matrices that yield even better bounds for the special cases $p = 1/3$ and $p = 2/3$.

\begin{theorem}\label{Thm:Upper}
For $p = \dfrac{a}{b}\leq \dfrac{1}{2}$, we have $\displaystyle V_2(p) \leq \frac{a}{b} - \left(1-\frac{a}{b}\right)^b \left(\frac{a}{b}\right)$.

For $p = \dfrac{a}{b}\geq \dfrac{1}{2}$, we have $\displaystyle V_2(p) \leq \frac{a}{b} - \left(1-\frac{a}{b}\right) \left(\frac{a}{b}\right)^b$.
\end{theorem}

The upper bounds from Theorem~\ref{Thm:Upper} for
selected rational values of $p$ are plotted as dots in Figure~\ref{Fig:Bounds}, with the lower bounds shown as a continuous curve obtained from Theorem~\ref{Thm:Distinct}.

\begin{figure}
\centering
\begin{tikzpicture}
\begin{axis}
\addplot[domain=0:0.5] {(x + x^2 + x^3 + 3*x^4 - 3*x^5 + x^6)/(2 + x + x^2 + x^3 - x^4)};
\addplot[domain=0.5:1] {(x + 5*x^2 - 10*x^3 + 10*x^4 - 5*x^5 + x^6)/(4 - 2*x - 2*x^2 + 3*x^3 - x^4)};
\addplot [only marks, mark size=0.5pt] table {
0.111111111111111 0.0726178426539313
0.125000000000000 0.0820488855242729
0.142857142857143 0.0942976175586980
0.166666666666667 0.110850337219936
0.200000000000000 0.134464000000000
0.250000000000000 0.170898437500000
0.333333333333333 0.234567901234568
0.500000000000000 0.375000000000000
0.666666666666667 0.567901234567901
0.750000000000000 0.670898437500000
0.800000000000000 0.734464000000000
0.833333333333333 0.777517003886603
0.857142857142857 0.808583331844412
0.875000000000000 0.832048885524273
0.888888888888889 0.850395620431709
};
\end{axis}
\end{tikzpicture}
\caption{Bounds on $V_2(p)$}
\label{Fig:Bounds}
\end{figure}

\section{Results for $n$ Players}\label{Sec:n_players}

\subsection{Terminology and Preliminary Results}\label{subsec:m_player_prelim}

The $n$-player Levine hat game
generalizes the 2-player game 
analyzed in the preceding sections.
In the general version of the game, each of the $n$ players can see all the hats on all of the other players' heads but not on his own. Once again, the players succeed if and only if each of them chooses a level on his own head that has a black hat (which 
we represent
by ``B'' or ``1,''
with a white hat being represented by ``W'' or ``0''). In the main version of this game, and throughout the current section unless stated otherwise, we suppose that the referee (or sultan, or warden) chooses all hat colors on all players' (or wise men's, or prisoners') heads independently to be black with probability $p = 1/2$. Within this section, $n$ will always refer to the number of players. Since it is sometimes important to distinguish clearly between one-person and multi-person subsets, 
we will think of the players as the wise men of Tanya Khovanova's problem statement~\cite{TK} and use singular masculine pronouns when referring to individual players.

As usual, for any fixed value of $n$, we suppose that there are $h$ hats per head ($h$ stands both 
for ``hats" and for ``height") and let $S$ be any joint strategy for the 
players.
Within this section 
we are usually interested
in $V_n(p)$, the optimal limiting $n$-player 
success probability (or ``value'')
as $h \rightarrow \infty$ for some given black-hat probability $p$, always with $0 < p < 1$. (As usual, 
we define
$q := 1-p$ throughout this section.) Often 
we are specifically interested
in $V_n$, the value of the game for the default case $p = 1/2$. As noted near the beginning of Section~\ref{Sec:2p}, other cases of particular interest have $p$ of the form $1/m$ (or $(m-1)/m$) for integer $m$, corresponding to a game in which each player must pick (or avoid) one of $m$ colors.

In the course of describing strategies or proving results, 
we sometimes refer
to restricted quantities such as
\[
V_n^{(h)}(p; S),
\]
the probability that the $n$ players win with a particular strategy $S$ when there are $h$ hats on each player's head.
It is not hard to see that $V_n^{(h)}(p)$ is nondecreasing in $h$ (we formally state this in Lemma~\ref{lem:nondecreasing_in_h}), so we 
have
\[
V_n(p) := \lim_{h \to \infty} V_n^{(h)}(p) = \sup_{h} V_n^{(h)}(p) .
\]

\begin{remark}
Recall from Section~\ref{sec:intro} that the optimal success probability, or value,
of any randomized strategy is always matched or exceeded by that of some deterministic
strategy. We therefore assume without loss of generality
that all players choose deterministic (pure) strategies, with each player's choice of level depending only on the (ordered) collection of $(n-1)h$ hat colors that he sees. Since, for each fixed $n$, $p$, and $h$, there are only finitely many possible strategies, the maximal value $V_ n^{(h)}(p)$ is actually achieved for some such strategy.
\end{remark}

\begin{lemma}\label{lem:nondecreasing_in_h}
We have
$V_n^{(1)}(p) = p^n$ for all $n \ge 1$, and the values $V_n^{(h)}(p)$ are nondecreasing in $h$ for each fixed $n$ and $p$.
\end{lemma}

\begin{proof}
The argument is straightforward. The players are free to ignore hats above any level, so increasing $h$ cannot hurt.
\end{proof}

\begin{lemma}\label{lem:V_n_nonincreasing}
For each $h \ge 1$ and $p \in (0,1)$, the values $V_n^{(h)}(p)$ and $V_n(p)$ are nonincreasing in $n$.
\end{lemma}

\begin{proof}
If an $(n+1)^{\rm st}$ player is added, his hat colors are independent of those for the first $n$ players and thus cannot help the first $n$ players to guess correctly with probability greater than $V_n^{(h)}(p)$ for finite $h$. Thus the full $(n+1)$-player set certainly cannot win with probability greater than $V_n^{(h)}(p)$. The result for $V_n(p)$ follows after taking suprema over $h \ge 1$.
\end{proof}

\subsection{Some General Bounds on $V_n(p)$}\label{subsec:n-player_gen_bounds}

One of our main results for the multiplayer game will be that $V_n := V_n(1/2)$ is actually {\em strictly} decreasing in $n$. Before proving this, however, 
we derive
some lower and upper bounds on $V_n(p)$, sometimes focusing on the case $p = 1/2$.

We already know from Lemma~\ref{lem:nondecreasing_in_h} 
that $V_n(p) \ge p^n$, but we will see now that we can do much better than random guessing of levels.

\begin{theorem}\label{thm:1_over_n+1}
We have
$V_n \ge 1/(n+1)$ and $V_n(p) \ge \frac{p/q}{n + p/q}$ for all $n \ge 1$.
\end{theorem}

\begin{proof}
Each player finds the first level at which the other $n-1$ players all have black hats, and he chooses that same level for his own head. (With probability 1 as $h \rightarrow \infty$, each player can find such a level.) The players succeed if and only if the first level that contains at least $n-1$ black hats actually contains $n$ black hats. There are $n+1$ possible arrangements of hats on this level, exactly one of them leads to success, and when $p = 1/2$ they are all equally likely, yielding the result for $p = 1/2$. For general $p$, the bound in the theorem follows from the posterior probability that the first candidate level actually has $n$ black hats.
\end{proof}

Thus we see that we can do dramatically better than with random guessing. In order to approach $1/(n+1)$ with the strategy above, however, we need $h$ to grow exponentially in $n$. The next strategy, discovered by Peter Winkler~\cite{TK} for the usual case $p = 1/2$ and generalized below for arbitrary $p$, requires $h$ to grow only logarithmically in $n$ and, quite surprisingly, yields another dramatic improvement in the success probability for large $n$.

\begin{definition}
For the next two theorems, given any $p \in (0,1)$, we let $q := 1-p$ as usual and define $r$ to be the reciprocal 
$r := 1/q$.
\end{definition}

\begin{theorem}\label{thm:winkler}
For every $p$ such that $0 < p < 1$, we have
\[
V_n(p) \ge \frac{1 - 1 / \ln(n)}{\lceil \log_{r} n + \log_{r} \ln n \rceil} \mbox{ for all } n \ge 3 .
\]
Thus for each $\epsilon > 0$, 
\[
V_n(p) \ge (1 - \epsilon) / \log_r(n) \mbox{ for all sufficiently large } n.
\]
\end{theorem}

\begin{proof}
Each player will attempt to choose the {\em first} level on which he has a black hat. (There might also be some serendipitous success probability coming from cases in which some players choose black hats but not the first black hats on their respective heads. We 
get a valid lower bound on $V_n$ by considering only the joint probability that every player chooses the {\em first} black hat on his own head.) We will further have each player choose only from the first $t$ levels, where $t = \lceil \log_{r} n + \log_{r} \ln n \rceil$. The players will hope that they each have at least one black hat within the first $t$ levels and that the sum of their $n$ first-black-hat levels is congruent (mod $t$) to some specified residue $s$. Given the values of $n$ and $t$, they will choose the residue $s$ during their strategy session to maximize their probability of hitting that residue. Conditioned on the assumption that all $n$ players have at least one black hat apiece within the first $t$ levels, the best residue will certainly occur with probability at least $1/t$.

The probability that any given player is bereft of black hats on the first $t$ levels is
\[
q^t \le q^{\log_{r} n + \log_{r} \ln n} = 1 / (n \ln n),
\]
and now, by a union bound, 
\[
\Prob(\mbox{At least one player has no usable black hats}) \le n / (n \ln n) = 1 / \ln(n).
\]
Thus $\Prob(\mbox{Every player has a usable black hat}) \ge 1 - 1 / \ln(n)$, which is positive if $n \ge 3$.

Conditioned on every player having a usable black hat, the best residue $s\pmod{t}$ for the sum of the $n$ lowest-black-hat levels occurs with probability at least $1/t$, and each player guesses the appropriate level on his own head to make the sum of all $n$ lowest-black-hat levels congruent to the target residue. The theorem follows. 
\end{proof}

In the next section 
we discuss several refinements to the basic Winkler strategy that improve the lower bounds by 5 to 10\% for moderate values of $n$ and outperform the $1/(n+1)$ lower bound for all $n \ge 3$. However, the following result (attributed in essence to Ori Gurel-Gurevich for the case $p = 1/2$ in a comment by hatmeister Lionel Levine on Tanya Khovanona's blog \cite{TK}, with no details given for the proof) shows that as long as each player is required to choose the {\em lowest} level at which he has a black hat, the lower bound from Theorem~\ref{thm:winkler} is asymptotically tight. Thus it seems likely that any substantial improvements upon the $1/ \log_r(n) = \ln(1/q) / \ln(n)$ approximate lower bound will require new ideas.

\begin{theorem}\label{thm:ori}
Suppose that each player must choose the {\em lowest} level on his own head that has a black hat (with the team failing if any player has only white hats). Then, using $\tilde{V}_n(p)$ to refer to the optimal probability as $h \rightarrow \infty$ that the players succeed (for any given black-hat probability $p \in (0,1)$) under this more stringent requirement, 
for each $\epsilon$ with $0 < \epsilon \le 1/4$ we have
\[
\tilde{V}_n(p) \le (1 + \epsilon) / \log_r(n) \mbox{ for all sufficiently large } n.
\]
\end{theorem}

\begin{proof}[Proof (Overview)]
We restrict attention to the first $t \approx \log_r(n)$ levels. An accomplice will uniformly at random choose a level $k$ from $\{1, \ldots, t\}$ and then uniformly at random choose a player $j$ from among those who happen to have lowest-black-hat level $Y_j$ equal to $k$. The accomplice will then inform the chosen player of his special status but will not tell him his lowest-black-hat level $k$. This special player is the only one required to guess his level $Y_j$; the accomplice will tell the other $(n-1)$ players (including those with $Y_j > t$) their own values $Y_j$.

As shown in the complete proof in Appendix~\ref{Subsec:ProofThmOri}, even with this substantial help (which the players can ignore if they wish), the one player who is not told his own level still cannot pick the correct level of his first black hat with probability much greater than $1/t$.
\end{proof}

For $n \ge 5$, the best lower bounds on $V_n$ that we know come from small refinements to Winkler's order-$(1 / \log(n))$ strategy. As $n$ grows, our best lower bounds on $V_n$ are asymptotically equal to $1 / \log_2(n)$. For upper bounds, we can use the fact that $V_n$ is nonincreasing in $n$ (from Lemma~\ref{lem:V_n_nonincreasing}) to see that
\[
V_n \le V_2 \le 81/224 \approx 0.361607 \mbox{ for all } n \ge 2,
\]
but we would like our upper bounds actually to (strictly) decrease as $n$ increases. A first step in this direction is the following theorem, which establishes a gap between $V_2$ and $V_3$.

\begin{theorem}\label{upperneq3}
We have
\[
V_3 \le 89/256 = 0.34765625 < 0.35 \le V_2 .
\]
\end{theorem}

\begin{proof}
We generalize the ``hint'' technique that was used earlier for 2 players to show that $V_2 \le 81/224.$ Given Players $A$, $B$, and $C$, we suppose that the referee gives $A$ a $2 \times 2$ hint with associated matrix that contains the two balanced binary columns of Hamming weight 1. In addition, the referee gives $B$ a $4 \times 6$ hint whose associated binary matrix contains all $\binom{4}{2} = 6$ distinct weight-2 columns of 4 bits each. The referee gives $C$ no hint.

It follows, much as for the 2-player upper bounds, that $C$ sees one of $6 \cdot 2 = 12$ column pairs for $(A, B)$ jointly and can then restrict attention to at most the first 12 levels on his own head, with each of the 12 possible $(A, B)$ column pairs associated with one of these levels on his own head. The (at most) 12 levels on $C$'s head can be permuted arbitrarily
without loss of generality,
leaving us with several million inequivalent strategies for $C$, which we put into the outer loop of a computer search. For each choice of $C$'s strategy, together with the realization of $C$'s hats within his (at most) 12 distinguished levels, $A$ and $B$ are left with a 2-player game with hints, and they choose their best (conditional) strategy so as to win with conditional probability $k / 12$ for some integer $k$.

Given $C$'s strategy and one of (at most) $2^{12}$ hat vectors on his distinguished levels, $B$ conditionally has only $4^2 = 16$ distinct strategies to consider, depending on which of the 4 rows he chooses on his own head for each of the 2 possible columns the referee could have placed on $A$'s head. For each joint choice of a $B, C$ strategy, together with what $A$ sees on $B$'s and $C$'s heads, $A$ guesses whichever of the 2 rows on his own head yields a higher (conditional) probability that $A, B,$ and $C$ will each point to a black hat (a `1'), breaking ties arbitrarily.

It turns out that, up to isomorphism, $C$'s unique best strategy uses only 6 levels on his own head. This best strategy can be represented by the following $2 \times 6$ matrix, indexed by the ``hint column"of $A$ and the hint column of $B$:

\[
\begin{pmatrix}
1 & 1 & 2 & 2 & 3 & 3\\
2 & 4 & 3 & 5 & 1 & 6
\end{pmatrix}.
\] 

When this optimal $C$-strategy is combined with optimal strategies for $B$ and $A$, the team wins with probability $267 / (64 \cdot 12) = 89/256 = 0.34765625$. Since we already know that $V_2 \ge 7/20 = 0.35$, we see that $V_2 - V_3 > 0.0023$, a nonzero gap.
\end{proof}

Shortly we will prove that $V_n$ is strictly decreasing in $n$. The proof will require a nontrivial upper bound on the probability that $n$ players can {\em avoid} all choosing black hats on their own heads, which one might call the ``mis{\`e}re" version of the game. (Stan Rabinowitz and Stan Wagon originally suggested this version of the 
problem.) More generally, one can try to bound 
\[
\Prob(\mbox{At least $k$ of the $n$ players choose black hats})
\]
or
\[
\Prob(\mbox{At most $k$ of the $n$ players choose black hats})
\]
for any value $k$ from 0 through $n$.

Still more generally, for any subset $\mathcal S$ of $\{0, 1, 2, \ldots, n\}$, one could ask for 
\[
\Prob(k \in \mathcal S), \mbox{ where $k$ is the number of players who choose black hats}.
\]
(E.g., $ \mathcal S$ might be the set of all even integers in $\{0, 1, 2, \ldots, n\}$.)

When $n = 2$, all of the problems above are equivalent to each other, as we now argue. 

For any value of $p$ (not necessarily $1/2$) and any joint strategy $S$ for the 2 players, 
we write
$\Prob(WW)$ for the probability that both players point to white hats, $\Prob(WB)$ for the probability that Players 1 and 2 point to white and black hats, respectively, on their own heads, etc. We use `*' as a wildcard; e.g., $\Prob(*B)$ is the probability that Player 1 points to a hat of either color and Player 2 points to a black hat.

Once we know $V_{2}(p; S)$ for any given strategy $S$, we know the entry $\Prob(BB)$ in a $2 \times 2$ matrix of probabilities with $\Prob(WW), \Prob(WB), \Prob(BW),$ and $\Prob(BB)$, and both row sums and both column sums are fixed (equal to $p$ for $\Prob(B*)$ and $\Prob(*B)$, equal to $1-p$ for $\Prob(W*)$ and $\Prob(*W)$). Thus the probability of every possible outcome or set of outcomes for any given strategy $S$ with 2 players can be calculated given $\Prob(BB)$. For example, our upper and lower bounds on $V_{2}(p)$ translate directly to upper and lower bounds on the probability that both players choose the same color hat, on the probability that the players choose hats of different colors, and on the probability that at least one player chooses a white hat (the latter being the mis{\`e}re version of the original game).

When $n > 2$, however, these different generalizations of the original game appear to be essentially different from one another. We will consider only the mis{\`e}re game (in which at least one player is supposed to choose a white hat), and only for the case $p = 1/2$. Letting $W_n$ be the maximal probability of winning the $n$-player mis{\`e}re game (where $W$ stands for ``white" and is the letter following the $V$ used for the original version of the game), we will see that $W_n \rightarrow 1$ as $n \rightarrow \infty$ but that $W_n$ is bounded away from 1 for each fixed $n$.

We begin by defining an infinite sequence of pairs of integers $(r_j, s_j)$ for $j \ge 1$ as below; these will be the dimensions of ``hint matrices" given to the various players for the $n$-player mis{\`e}re game.

\begin{definition}
Let $(r_1, s_1) = (2,2), \ (r_2, s_2) = (4, 6)$, and for
each $k \ge 3$, define 
\[
r_k = 2 \prod_{j=1}^{k-1} s_j \mbox{ and } s_k = \binom{r_k}{r_k / 2}.
\]
\end{definition}

\begin{theorem}\label{thm:misere}
Letting $W_n = \sup_S \Prob(\mbox{At least 1 of $n$ players chooses a white hat})$ over all possible $n$-player strategies as $h$, the number of hats per head, grows to infinity, we have the following:
\[
W_2 = 1/2 + V_2 \in [17/20, 193/224] \mbox{ (hence in [0.85, 0.8616...])}
\]
and, for $n \ge 3$,
\[
1 - (1/2)^n \le W_n \le 1 - \frac{1}{(r_{n} / 2) 2^{(r_{n} / 2)} }.
\]
\end{theorem}

\begin{proof}
It is straightforward to show that 
\[
W_2 = \sup_S ({\Prob(WW) + \Prob(WB) + \Prob(BW)}) = \sup_S ({\Prob(W*) + \Prob(BW)}).
\]
For $p = 1/2$, though, $\Prob(W*) = 1/2$ for every strategy $S$, and since $p = 1-p$,
\[
\sup_S \Prob({BW}) = \sup_S \Prob({BB}) = V_2,
\]
and the bounds for $W_2$ follow immediately.

The lower bounds on $W_n$ for all $n \ge 3$ follow by letting each player just choose the first level on his own head. The upper bounds follow from a cascading ``hint" technique generalizing the upper-bounding techniques used earlier for $V_2$ and $V_3$. The remaining details of the proof are 
given in
Appendix~\ref{app:subsec:W_n_lt_1}.
\end{proof}

Now that we have bounded $W_n$ away from 1 for every $n$ (for the mis{\`e}re game), we are finally ready to show that $V_n$ (for the original game) is strictly decreasing in $n$.

\begin{theorem}\label{thm:strictly_decreasing}
We have 
$V_{n+1} < V_n$ for all $n \ge 1$.
\end{theorem}

\begin{proof}
Given any $h$ and strategy $S$ for a game with $n+1$ players, let us write $\Prob(B \ldots B,B)$ to refer to the probability that Players 1 through $n+1$ all choose black hats. Let us write $\Prob(B \ldots B,W)$ for the probability that the first $n$ players choose black hats while Player $n+1$ chooses a white hat, and write
\[
\Prob(B \ldots B,*) := \Prob(B \ldots B,B) + \Prob(B \ldots B,W)
\]
for the probability that the first $n$ players choose black hats while Player $n+1$ chooses a hat of arbitrary color.

Now, for any given $(n+1)$-player strategy, we have
\[
\Prob(B \ldots B,B) = \Prob(B \ldots B,*) - \Prob(B \ldots B,W).
\]

Thus, taking suprema over strategies and over $h$, we have
\begin{align*}
V_{n+1} & = \sup \Prob(B \ldots B,B) \\ 
& = \sup \ \{ \Prob(B...B,*) - \Prob(B...B,W) \} \\
& \le \sup \Prob(B \ldots B,*) + \sup \ \{1 - \Prob(B \ldots B,W) \} - 1 \\
& = V_n + \sup \ \{1 - \Prob(B \ldots B,B) \} - 1,
\end{align*}
where the last equality is because $p = 1/2$, so the maximal probability of avoiding any given ordered sequence of chosen hat colors for the $n+1$ players is the same as the maximal probability for any other ordered sequence of chosen hat colors.

Now we have
\[
V_{n+1} \le V_n + W_{n+1} - 1 = V_n - (1 - W_{n+1}),
\]
so $(V_n - V_{n+1}) \ge (1 - W_{n+1})$, which by Theorem~\ref{thm:misere} is 
positive and bounded away from 0 for each fixed $n$.
\end{proof}

\begin{remark}
A similar argument shows that $W_n$ for the mis{\`e}re game is strictly increasing in $n$. See Theorem~\ref{thm:misere_increasing} in Appendix~\ref{App:n_player_results} for the proof.
\end{remark}

One of our main unresolved questions 
concerns Levine's original conjecture, which we restate here.

\begin{repconjecture}{Conjecture:Levine}[Levine]
The optimal success probability in the $n$-player game is $o(1)$
as $n$ goes to infinity; i.e.,
\[ \lim_{n \to \infty} V_n = 0. \]
\end{repconjecture}

\begin{remark}
Although the conjecture above seems likely to be true, the rate of decrease provided in Theorem~\ref{thm:strictly_decreasing} is insufficient to prove it, even with tighter bounds on $(1 - W_n)$. The upper bound that Theorem~\ref{thm:strictly_decreasing} yields on $\lim\limits_{n\to\infty} V_n$ is
\[
V_2 - \sum_{k=3}^\infty (1-W_k),
\]
which is at least
\[
\frac{7}{20} - \sum_{k=3}^\infty \frac{1}{2^k} = \frac{7}{20} - \frac{1}{4} = \frac{1}{10} .
\]
\end{remark}

\section{Best Current Lower Bounds}\label{Sec:best_lower_bounds}
For puzzle aficionados (or prisoners' advocates) who would like to improve upon existing strategies for various numbers of players (or prisoners), we collect here the best lower bounds we know on $V_2$ through $V_{12}$ (all for $p = 1/2)$. The result for $V_2$ comes from Section~\ref{Sec:2p} and is conjectured to be optimal. The results for $V_3$ and $V_4$ essentially come from hill climbs over symmetric strategies with 5 or 4 hats, respectively, per player.

The bounds on $V_5$ through $V_{12}$ come from generalized Winkler-style strategies (as described in Theorem~\ref{thm:winkler}) in which the players focus on the first $t \approx \log_2(n)$ levels and all try to select the {\em lowest} levels on which they have black hats. However, the $n$ players use general $t \times t$ Latin-square operations, not necessarily mod-$t$ addition, in order to construct a ``sum" in $\{0, 1, \ldots, t-1\}$ of their $n$ respective lowest-black-hat levels. Furthermore, in assessing each candidate strategy, we take into account all ``bonus" success probability that arises when the players miss their target ``sum" (mod $t$) but nonetheless all serendipitously point to black hats.

Finally, all $t$-level strategies --- whether found by hill climbing or derived from Latin-square operations --- are then augmented by working with $t$-level ``tiers" and having each player recursively ``reset" (shifting up $t$ levels at a time) whenever he sees at least one other player with an all-white stack of hats in the current tier. A second, smaller, improvement comes from also recursively resetting whenever a player sees {\em only} black hats on {\em all} other players' heads within the current tier. Using both ``white" and ``black" recursive resets leads to rational lower bounds with denominators of the form
\[
(2^t - 1)^n + n (2^t - 1)^{n-1} - (n+1),
\]
and we retain these unreduced fractions in the table below.

Additional details about the best strategies known are given in 
Appendix~\ref{App:subsec:details_best_strategies}.
The resulting values of $V_n$ for $2 \le n \le 12$ are shown 
below, with all decimal values rounded down to 6 decimal places to provide true lower bounds:

\[
\begin{array}{lcccl}
V_2 & \ge & 21/60 & = & 0.350000 , \\
V_3 & \ge & 9119/32670 & = & 0.279124... , \\
V_4 & \ge & 14844/64120 & = & 0.231503... , \\
V_5 & \ge & 205447/1012494 & = & 0.202911... , \\
V_6 & \ge & 2984604 / 15946868 & = & 0.187159... , \\
V_7 & \ge & 43930663 / 250593742 & = & 0.175306... , \\
V_8 & \ge & 651583632 / 3929765616 & = & 0.165807... , \\
V_9 & \ge\ & & & 0.158764... , \\
V_{10} & \ge & & & 0.153517... , \\
V_{11} & \ge & & & 0.149025... , \\
V_{12} & \ge & & & 0.145047... . \\
\end{array}
\]

\section{Future Directions}

This paper leaves open certain questions that might be of interest to other 
hatters, mad or otherwise.
In particular, there are two conjectures that
we stated earlier, phrased below as questions:

\begin{enumerate}
\item Is $V_2$ exactly equal to 0.35?

\item Does $V_n$ approach 0 as $n\to\infty$?
\end{enumerate}

We also note that for the upper-bound results in Section~\ref{Sec:Matrix}, one promising family of matrices has size $2^k\times (2^{k+1}-2)$; the cases for $k=1,2,$ and $3$ are presented in that section. For any such matrix, one must consider $\mathrm{Bell}(2^{k+1}-2)$ partitions of the columns, a task that seems infeasible beyond $k=3$. One might look for ways to reduce the computational difficulty of this approach. For example, we can safely assume that none of the column subsets within an optimal partition contains both a column and its bitwise complement. Unfortunately, the resulting reduction in work is not very significant.

There are also several generalizations 
one could consider,
some of which are presented below:

\begin{enumerate}
\item In the $n$-player case, one might require at least $k$ players to pick a black hat (the $k=n$ and $k=1$ cases are discussed above).

\item 
One could allow more than 2 hat colors,
perhaps with a different payoff system.
For example, one might consider a grayscale version of the game, where each hat 
has
a value in $[0,1]$, with $0$ being white and $1$ being black, and with the value of a joint guess taken to be the product of the values of the selected hats.
\end{enumerate}

\bigskip
\bigskip

\noindent {\bf Acknowledgements.}
We have corresponded or spoken with Aaron Atlee, Larry Carter, Joseph DeVincentis, Eric Egge, 
Ehud Friedgut, Jerry Grossman, Gil Kalai, Tanya Khovanova, Sandy Kutin, Lionel Levine,
Stephen Morris, Rob Pratt, Jay-C Reyes, Joel Rosenberg, Walter Stromquist, Alan
Taylor, Dan Velleman, Stan Wagon, Peter Winkler, Chen Yan, Piotr Zielinski, and
no doubt others. 
The eighth author received funding from Excellence of Science grant number 30950721, ``Symplectic Techniques."

\newpage

\appendix

\section{The 1-Player Game}\label{App:pseudo-theorem}

Now we consider the curious case of the 1-player (or ``solo'')
version of the Levine hat game,
in which Sol must try to point to a black hat on his own head.  
Imagine, if you will, that Sol can’t see any of the infinitely many black and 
white hats on his own head, but that he is endowed with infinite 
computational abilities and armed with a secure faith in the Axiom of Choice ({\bf AC}).  
In order to describe Sol's strategy, 
we first look
at an auxiliary game with countably many players.

We begin with some definitions. Throughout this section, 
$I$ denotes
a fixed countable set, which we call the {\em set of players}. A {\em hat-stack} is an infinite sequence of zeros and ones, or if preferred, an infinite sequence from $\{ \text{\it white}, \text{\it black}\}$.  Each 
member
of a hat-stack will be called a {\em hat}. A {\em hat-assignment} is a mapping that assigns a hat-stack to each element of $I$. Note that the image of a hat-assignment is a set of hat-stacks. Given a set of hat-stacks, $M$, a {\em black level} of $M$ is a natural number, $i$, such that $m(i) = 1$ for every $m \in  M$. In addition, $M$ will be called {\em generic} if every finite subset of $M$ has a black level, and $M$ will be called {\em almost-generic} if some cofinite subset of $M$ is generic.

In a hat-assignment there are countably many hats (each element of $I$ is assigned a hat-stack, and each hat-stack has countably many hats). A hat-assignment gives each of these hats a color. A cylinder is an assignment of colors to some finite subset of these hats. The standard Bernoulli measure ($p = 1/2$) is a probability measure defined on the $\sigma$-algebra of sets generated by the cylinders. 
With this probability space, the following proposition is standard and is stated without proof.

\begin{proposition} In the standard Bernoulli measure ($p = 1/2$), almost every hat-assignment is one-to-one 
(i.e., with no two elements of $I$ receiving the same hat-stack) and has a generic image.
\end{proposition}

\begin{theorem} Assume the Axiom of Choice. Then there exists a set $C$ whose elements are countable sets of hat-stacks, and $C$ has the following property: For every countable set of hat-stacks, $T$, there is a unique element $R \in C$ such that
\begin{enumerate}
    \item $T \backslash R$ is finite, and
    \item $R \backslash T$ is finite.
\end{enumerate}
Furthermore, if $T$ is generic, then $R \backslash T$ has a black level.
\end{theorem}

\begin{proof} Let two countable sets of hat-stacks be equivalent if they differ by a finite set of hat-stacks. Using {\bf AC}, let $E$ be a choice set for the set of equivalence classes. Using {\bf AC} a second time, replace each $D \in E$ that is almost-generic, with a cofinite generic subset of $D$. Since this does not change the equivalence class of $D$, the resulting set $C$ is also a choice set. In addition, any 
almost-generic
element of $C$ is generic. Therefore, given a countable set of hat-stacks, $T$, there is a unique element $R \in C$ that is equivalent to $T$. This gives the first two properties. Furthermore, if $T$ is generic, then any subset of $T$ is also generic. So, removing the finite set $R \backslash T$ from $R$ results in the generic set $R \cap T$. Therefore, $R$ is 
almost-generic.
But $R \in  C$, so $R$ is generic. Since $R \backslash T$ is a finite subset of $R$, it has a black level.
\end{proof} 

We now interpret the previous results as a hat game.

The {\em Auxiliary Game} has a countably infinite set of players, each with an infinite sequence of hats. The usual rules apply, including ``No looking at your own hats.'' The players win if all but finitely many of them are able to point to a black hat.  This game is 
similar 
to a puzzle described by Greg Muller \cite{GM}, who attributed the earlier puzzle to Mike O'Connor.
In that game there is a fixed ordering of the players.  Here we want the players to be indistinguishable.

Our previous results
provide a strategy for winning the Auxiliary Game with probability 1.
Let $I$ be the set of players. Consider a random hat-assignment, chosen according to the probability space defined above. Using the proposition and ignoring a measure-zero event, 
we find that
the assignment is one-to-one and its image, $T$, is generic. 
Let $C$ be as in the previous theorem, and let $R \in C$ be the unique set satisfying the conclusions of the theorem.

To describe the strategy, fix a player $i \in I$. Player $i$ cannot determine $T$, since 
no players see
their own hats. Nevertheless, Player $i$ can determine $T_i$, 
the set of hat-stacks 
assigned to the other players.
The theorem applies to $T$ as well as to $T_i$, and by uniqueness, both yield the same element $R$. So, Player $i$ reports a black level of the finite set $R \backslash T_i$.

Suppose that each player follows the strategy above. If Player $i$ has a hat-stack in $R$, then since the assignment is one-to-one, this hat stack will also be in $R \backslash T_i$. By choosing a black level of $R \backslash T_i$, 
any such player is guaranteed to choose a black hat.
Since $T \backslash R$ is finite, the game is won.

The Auxiliary Game also shows how Sol can win the solo game. On the fateful day, he brings with him infinitely many friends, each with his own referee and coin. This is not explicitly against the rules. 
All referees simultaneously select 
random hat sequences for their respective players. Sol and his friends play the Auxiliary Game. With probability 1, all but finitely many players choose a black hat. Sol is confident that he will not be one of the unlucky ones. Now, he might have a persnickety logician friend who warns against depending on an event whose probability cannot be precisely measured.  Sol should ignore this advice! Since only finitely many players fail, and all players are essentially identical, he is virtually guaranteed to win the game!

\section{Computing the Performance of 2-Player Strategy $S_1$}\label{Sec:S1App}

Below we compute the success probability, or value, $V_2(p; S_1)$ for the strategy $S_1$ from Section~\ref{Sec:2p}.
We break the calculation of $V_2(p; S_1)$ into 
7 cases.

Case 1: Both players are mono\-chromatic to the same odd position. By ``mono\-chromatic up to an odd position $2\ell +1$,'' we mean that a player has either all W or all B up to position $2\ell +1$ but not up to position $2\ell+3$.

Case 1(a): Both players start with B and have B hats up to position $2\ell+1$ but not to position $2\ell+3$, for some integer $\ell$. The probability of winning conditioned on this is given 
in Table~\ref{Table:1a}, which shows hats in positions $2\ell+1$, $2\ell+2$, and $2\ell+3$ for both players.

\begin{table}[hb]
\centering
\begin{tabular}{|c||c|c|c|}
\hline
 & $\blacksquare\blacksquare\square$ & $\blacksquare\square\blacksquare$ & $\blacksquare\square\square$ \\
\hline\hline
$\blacksquare\blacksquare\square$ & $p^2q^2$ & $p^2q^2$ & $pq^3$ \\
\hline
$\blacksquare\square\blacksquare$ &$p^2q^2$ & & \\
\hline
$\blacksquare\square\square$ & $pq^3$ & & $q^4$ \\
\hline
\end{tabular}
\caption{Case 1(a)}\label{Table:1a}
\end{table}

The probability of this case occurring and 
the players' winning is
\[
p^2(3p^2q^2 + 2pq^3 + q^4) + p^6(3p^2q^2 + 2pq^3 + q^4) + p^{10}(3p^2q^2 + 2pq^3 + q^4) + \dots ,
\]
which we can simplify by summing the geometric series to obtain
\[
\frac{(pq)^2(1 + 2p^2)}{1 - p^4}.
\]

Case 1(b): Both players start with W and have W hats up to position $2\ell+1$ but not to position $2\ell+3$, for some integer $\ell$. The probability of winning conditioned on this is
given in Table~\ref{Table:1b}, which shows hats in positions $2\ell+1$, $2\ell+2$, and $2\ell+3$.

\begin{table}[hb]
\centering
\begin{tabular}{|c||c|c|c|}
\hline
 & $\square\blacksquare\blacksquare$ & $\square\blacksquare\square$ & $\square\square\blacksquare$ \\
\hline\hline
$\square\blacksquare\blacksquare$ & $p^4$ & & $p^3q$ \\
\hline
$\square\blacksquare\square$ & & & $p^2q^2$ \\
\hline
$\square\square\blacksquare$ & $p^3q$ & $p^2q^2$ & \\
\hline
\end{tabular}

\caption{Case 1(b)}\label{Table:1b}
\end{table}

The probability of this case occurring and 
the players' winning is
\[
q^2(p^4 + 2p^3q + 2p^2q^2) + q^6(p^4 + 2p^3q + 2p^2q^2) + q^{10}(p^4 + 2p^3q + 2p^2q^2) + \dots ,
\]
which simplifies to
\[
\frac{(pq)^2}{1 - q^2}.
\]

Case 1(c): One player starts W and the other starts B, or vice versa, and both are 
monochromatic
to the same odd position.

\begin{table}[ht]
\centering
\begin{tabular}{|c||c|c|c|}
\hline
 & $\blacksquare\blacksquare\square$ & $\blacksquare\square\blacksquare$ & $\blacksquare\square\square$ \\
\hline\hline
$\square\blacksquare\blacksquare$ & & $p^3q$ & \\
\hline
$\square\blacksquare\square$ & & $p^2q^2$ & \\
\hline
$\square\square\blacksquare$ & & & \\
\hline
\end{tabular}

\caption{Case 1(c)}\label{Table:1c}
\end{table}

The probability of this case occurring and 
the players' winning 
is
\[
pq(p^3q + p^2q^2) + p^3q^3(p^3q + p^2q^2) + p^5q^5(p^3q + p^2q^2) + \dots ,
\]
which simplifies to
\[
\frac{(pq)^2(p^2 + pq)}{1 - (pq)^2} = \frac{p^3q^2}{1 - (pq)^2}.
\]
Finally, notice that the roles of the two players could be interchanged here, so we double the above probability to get
\[
\frac{2p^3q^2}{1 - (pq)^2}.
\]

Case 2: Players are monochromatic to different odd positions.

Case 2(a): The taller monochromatic stack is W. In this case the player with the taller W stack will always choose a W hat, so the probability of winning is 0.

Case 2(b): The taller monochromatic stack is B.
 
Case 2(b)(i): The taller monochromatic B stack is taller than the shorter one by at least 2 odd positions. In this case, the player with the taller B stack always guesses correctly. The player with the shorter stack guesses correctly with probability $p$ since his guess is uncorrelated with the other player's guess. To calculate the probability of winning in this case, note that the probability of a player being monochromatic up to odd position $2\ell+1$ and not to $2\ell+3$ is $p^{2\ell+1}(1 - p^2) + q^{2\ell+1}(1 - q^2)$. The probability of the other player being monochromatic B up to at least position $2\ell+5$ is $p^{2\ell+5}$. Thus, the probability of this case occurring and 
the players' winning is
\[
\sum_{\ell = 0}^\infty (p^{2\ell+1}(1 - p^2) + q^{2\ell+1}(1 - q^2))p^{2\ell+5},
\]
which simplifies to
\[
\frac{p^6(1 - p^2)}{1 - p^4} + \frac{p^5q(1 - q^2)}{1 - (pq)^2}.
\]
Thus, the probability of winning in this case is obtained by multiplying by $p$ and by 2, giving
\[
\frac{2p^7(1 - p^2)}{1 - p^4} + \frac{2p^6q(1 - q^2)}{1 - (pq)^2}.
\]

Case 2(b)(ii): The taller stack is 
monochromatic
B to odd position $2\ell+3$ and the shorter stack is 
monochromatic
B to position $2\ell+1$, for some integer $\ell$. The player with the taller B stack always guesses correctly. The player with the shorter B stack guesses correctly according to Table~\ref{Table:2bii}, which shows hats in position $2\ell+1$ to $2\ell+5$ for the player with the taller B stack, and $2\ell+1$ to $2\ell+3$ for the other player. Note that in the middle row, the player with the shorter stack guesses the hat in position $2\ell + 4$ and so has probability $p$ of being correct, which is multiplied by the probability of the situation occurring.

\begin{table}[hb]
\centering
\begin{tabular}{|c||c|c|c|}
\hline
 & $\blacksquare\blacksquare\square$ & $\blacksquare\square\blacksquare$ & $\blacksquare\square\square$ \\
\hline\hline
$\blacksquare\blacksquare\blacksquare\blacksquare\square$ & & $p^2q^2$ & \\
\hline
$\blacksquare\blacksquare\blacksquare\square\blacksquare$ & $p^3q^2$ & $p^3q^2$ & $p^2q^3$ \\
\hline
$\blacksquare\blacksquare\blacksquare\square\square$ & & $pq^3$ & \\
\hline
\end{tabular}
\caption{Case 2(b)(ii)}\label{Table:2bii}
\end{table}

If this case occurs, then the probability of winning is given by
\[
p^2q^2 + 2p^3q^2 + p^2q^3 + pq^3 = pq^2(1 + p + p^2).
\]
The probability of this case occurring and 
the players' winning
is, therefore, given by
\[
p^4pq^2(1 + p + p^2) + p^8pq^2(1 + p + p^2) + p^{12}pq^2(1 + p + p^2) + \dots ,
\]
which simplifies to:
\[
\frac{p^5q(1 - p^3)}{1 - p^4}.
\]
Taking into account 
the fact
that either player could have the taller stack, we get a probability of winning as
\[
\frac{2p^5q(1 - p^3)}{1 - p^4}.
\]

Case 2(b)(iii): The taller B stack is 
monochromatic
B to odd position $2\ell+3$ and the shorter stack is 
monochromatic
W to position $2\ell+1$, for some integer $\ell$. The player with the taller B stack always guesses correctly. The player with the shorter W stack guesses correctly according to Table~\ref{Table:2biii}.

\begin{table}[ht]
\centering
\begin{tabular}{|c||c|c|c|}
\hline
 & $\square\blacksquare\blacksquare$ & $\square\blacksquare\square$ & $\square\square\blacksquare$ \\
\hline\hline
$\blacksquare\blacksquare\blacksquare\blacksquare\square$ & $p^3q$ & & $p^2q^2$ \\
\hline
$\blacksquare\blacksquare\blacksquare\square\blacksquare$ & $p^4q$ & $p^3q^2$ & $p^3q^2$ \\
\hline
$\blacksquare\blacksquare\blacksquare\square\square$ & $p^2q^2$ & & $pq^3$ \\
\hline
\end{tabular}

\caption{Case 2(b)(iii)}\label{Table:2biii}
\end{table}

If this case occurs, then the probability of winning is given by
\[
p^4q + 2p^3q^2 + p^3q + 2p^2q^2 + pq^3 = pq(1 + p - pq^2).
\]
The probability of this case occurring and 
the players' winning
is, therefore, given by
\[
p^3qpq(1 + p - pq^2) + p^5q^3pq(1 + p - pq^2) + p^7q^5pq(1 + p - pq^2) + \dots ,
\]
which simplifies to
\[
\frac{p^4q^2(1 + p - pq^2)}{1 - (pq)^2}.
\]
Taking into account 
the fact
that either player could have the taller stack, we get a probability of winning as
\[
\frac{2p^4q^2(1 + p - pq^2)}{1 - (pq)^2}.
\]

Summing the success probabilities from all the cases above and replacing $q$ with $1 - p$ yields
\[
V_2(p;S_1) = \frac{p (1 + p + p^2 + 3 p^3 - 3 p^4 + p^5)}{2 + p + p^2 + p^3 - p^4}.
\]

\section{Proof That the Matrix-Based Upper and Lower Bounds Converge to $V_2(p)$}\label{Sec:MatrixBounds}

For each rational probability $p=a/b$ in lowest terms and each positive integer $m$ that is a multiple of $b$, we define two matrices, $L_{m,p}$ and $U_{m,p}$. The columns of these matrices will be elements of $\{0,1\}^m$. In $L_{m,p}$ all $2^m$ such columns appear, and each column with $t$ 
1s occurs $a^t (b-a)^{m-t}$ times, for a total of $b^m$ columns. In $U_{m,p}$ only the $\binom{m}{mp}$ columns with $mp$ 
1s
occur, and there is no repetition of columns. For both matrices, the columns may be ordered arbitrarily.

As argued in Section~\ref{Sec:Matrix}, $V(L_{m,p})$ is a lower bound for the value of the two-person hat game with black-hat probability $p$, and $V(U_{m,p})$ is an upper bound. Below 
we prove
Theorem~\ref{Thm:UpperLower}, the main theorem of that section, that as $m \rightarrow \infty$, the matrix-based upper and lower bounds converge to $V_2(p)$, the value of the 2-player hat game with (rational) black-hat probability $p \in (0,1)$.

The following three lemmas establish background results needed to prove Theorem~\ref{Thm:UpperLower}.

\begin{lemma}[Determinism Lemma]\label{Lemma:Determinism}
If there are duplicate columns in $M$, then Player 2 may as well put them into the same equivalence class when playing the matrix game on $M$.
\end{lemma}
	
\begin{proof}
This is just a simple convexity argument. More precisely, let $\sim$ be an equivalence relation chosen by Player 2. Let $c$ and $d$ be two duplicate columns in two different equivalence classes, $C$ and $D$, resp. Consider the two equivalence relations, $\sim_C$ and $\sim_D$ which are identical to $\sim$ except that in $\sim_C$ we move column $d$ to $C$ and in $\sim_D$ we move column $c$ to $D$. Let $v$ be any vector that respects $\sim$. We form two new vectors $v^C$ and $v^D$ that are identical to $v$ except in positions $c$ and $d$, where we have $v_c^C=v_d^C=v_c$ and $v_c^D=v_d^D=v_d$. Note that $v^C$ respects $\sim_C$ and $v^D$ respects $\sim_D$. Note also that these mappings might not be one-to-one; if $d$ is the only member of $D$, then it is possible that $v^C = w^C$ even though $v \neq w$. In this case we will count $v^C$ and $w^C$ as two different vectors, so we can write $\Prob(v) = \Prob(v^C) = \Prob(v^D)$. Let $r$ be any row of $M$. Then $r \cdot c = r \cdot d$ because $c$ and $d$ are duplicate columns. 
	
So $r \cdot v = (r \cdot v^C + r \cdot v^D)/2$. If we let $r(v)$ denote the row that Player 1 assigns to $v$, then we have
\[
2 \sum\limits_{v} (r(v) \cdot v) \Prob(v) = \sum\limits_{v} (r(v) \cdot v^C) \Prob(v^C) + \sum\limits_{v} (r(v) \cdot v^D) \Prob(v^D).
\]

Therefore, at least one of sums on the right is at least as big as the (undoubled) sum on the left. Assume, without loss of generality, that 
\[
\sum_{v} (r(v) \cdot v) \Prob(v) \leq \sum_{v} (r(v) \cdot v^C) \Prob(v^C).
\]
Then $V(M; \sim) \leq V(M; \sim^C )$.

\end{proof}

\begin{lemma}[Replication Lemma] \label{Lemma:Duplication}
If all columns are replicated the same number of times, this does not change the value of the matrix.
\end{lemma}

\begin{proof}
Suppose $M_2$ is an $m \times kn$ matrix formed from the $m \times n$ matrix $M_1$ by including each column of $M_1$ a total of $k$ times. Let $\sim_2$ be an equivalence relation chosen by Player 2 on $M_2$. By the Determinism Lemma we can assume that $\sim_2$ assigns each set of identical columns to the same equivalence class. Let $\sim_1$ be the restriction of $\sim_1$ to the original matrix $M_1$. Let $v_1$ be any vector of numbers that respects $\sim_1$, and let $v_2$ be the $k$-fold expansion of $v_2$. Let $r$ be a row of $M_2$. Then $r \cdot v_2 = kr \cdot v_1$ Therefore, $V(M_2;\sim_2) = V(M_1; \sim_1)$.
\end{proof}

\begin{lemma}[Erasure Lemma] \label{Lemma:Erasure}
If we remove a small proportion $\epsilon$ of the columns of $M$, the value $V(M)$ changes by at most $\frac{2 \epsilon}{1-\epsilon}$.
\end{lemma}

\begin{proof}
Let $r$ be any row of $M$ with $n$ columns, and $v$ any vector of size $n$ whose entries are zeros and ones. Then when the columns are removed, $r \cdot v$, which is at most $n$, becomes $r' \cdot v'$ and decreases by some amount $k$, where $0\leq k \leq \epsilon n$. So
\begin{align*}
\left| \frac{r \cdot v}{n} -\frac{r' \cdot v'}{n(1-\epsilon)} \right|&= \left|\frac{r \cdot v}{n} -\frac{r \cdot v - k}{n(1-\epsilon)} \right| \\
&= \left| \frac{ -\epsilon r \cdot v +k}{n(1-\epsilon)} \right| \\
&\leq \left| \frac{- \epsilon}{1- \epsilon} \right|+ \left|\frac{\epsilon}{1- \epsilon} \right| \\
&= \frac{2 \epsilon}{1- \epsilon} .
\end{align*}
\end{proof}

\begin{lemma}[Perturbation Lemma] \label{Lemma:Perturbation}
Let $\epsilon>0$ and let $M$ and $N$ be two matrices of the same size with elements in $\{0,1\}$. Let 
\[
E := N-M = 
\begin{pmatrix}
\epsilon_{1,1} & \epsilon_{1,2} & \cdots & \epsilon_{1,n} \\
\epsilon_{2,1} & \epsilon_{2,2} & \cdots & \epsilon_{2,n} \\
\vdots & \vdots & \ddots & \vdots \\
\epsilon_{m,1} & \epsilon_{m,2} & \cdots & \epsilon_{m,n} 
\end{pmatrix},
\]
with each entry $\epsilon_{i,j} \in \{-1,0,1\}$. If the average of the absolute values of the entries in each row of $E$
satisfies
\[
\frac{1}{n} \sum\limits_{j=1}^{n} |\epsilon_{i,j} | < \epsilon,
\]
then $|V(M + E) - V(M)| < \epsilon$.
\end{lemma} 

\begin{proof}
When $E$ is added to $M$, each dot product $v \cdot r_i$ changes by less than $n \epsilon$. Thus the value of the best dot product for $v_i$ changes by less than $n \epsilon$. So the value of the matrix changes by less than $\frac{n \epsilon}{n}=\epsilon$.
\end{proof}

We now prove that the upper and lower matrix-based bounds converge for the 2-player game.

\begin{reptheorem}{Thm:UpperLower}[Convergence Theorem]
Let $L_{m,p}$ and $U_{m,p}$ be defined as above. Then
\[
\lim_{m \rightarrow \infty} V(U_{m,p}) - V(L_{m,p}) = 0.
\]
\end{reptheorem}

\begin{proof}
We first give a brief outline. We shall begin with $L_{m,p}$ and, using the Replication, Erasure, and Perturbation lemmas, will move to $U_{m,p}$ and notice that the value of $V$ will not have changed much. Each entry of $L_{m,p}$ is in $\{0, 1\}$. We will first remove a small proportion of the columns from $L_{m,p}$, ones that are far out of balance. Then we will replicate the remaining columns (each column replicated the same number of times). Then we will flip some values of the remaining columns, in order to bring them into balance. Each entry is changed by a perturbation $\epsilon_{i,j} \in \{ 1, 0, -1\}$. In order to appeal to the Perturbation Lemma, we will need to keep the average absolute values of these perturbations small along each row, which is the same as keeping the proportion of changes small along each row.
	
But we will keep a symmetry of the rows, so that each row will receive the same number of changes. So it suffices to keep small the proportion of the matrix that is changed. The resulting matrix will be a replication of $U_{m,p}$, finishing the proof.
	
More precisely, let us call a column {\em balanced} if the proportion of ones in the column is $p$. Fix $\epsilon$ and using the Law of Large Numbers, let $m$ be large enough that a proportion of at most $\epsilon$ of the columns of $L_{m,p}$ are not within $m \epsilon$ bits being balanced (i.e., if 1 appears in a column $t$ times then $|t - mp| \geq m \epsilon$). We delete any column that satisfies this condition. This deletes at most a proportion $\epsilon$ of the columns and so	by the Erasure Lemma, this changes the value $V$ by at most $2\epsilon/(1-\epsilon)$. For each of the remaining columns, replicate it a total of $m!$ times. By the Replication Lemma, this does not change the value of the matrix. Now replace each of these duplicates by one of its nearest balanced neighbors. For a column with $t$ `1' bits the number of balanced neighbors to choose from is $\binom{t}{t-mp}$ in the case that $t \geq mp$, and $\binom{m-t}{mp-t}$ otherwise. In any case, the number of choices divides $m!$. Therefore, we can make sure that each neighbor is used the same whole number of times. 

Recall that in $L_{m,p}$, each column with the same number of ones occurred the same number of times. By construction, the same will be true of our resulting matrix, $N$. In other words, $N$ is just a replication of $U_{m,p}$. Furthermore, the maximum number of changes made to any column is bounded by $m\epsilon$. By symmetry, each row will receive the same number of changes, so the proportion of changes in each row is also at most $\epsilon$. Therefore, the Perturbation Lemma applies, and $|V(U_{m,p}) - V(L_{m,p})| \leq 2\epsilon/(1-\epsilon) + \epsilon$. Since $\epsilon$ is arbitrary, 
the theorem follows.
\end{proof}

\section{Upper Bounds on $V_2(p)$ for Rational $p$}\label{App:SpecificUpperBounds}

In this appendix we discuss some specific upper-bound results for the 2-player game, along with proving the upper bounds claimed in Theorem~\ref{Thm:Upper}.

\subsection{Dual Strategies}\label{SubSec:Dual}

One tool that is useful for both upper and lower bounds 
is the notion 
of a dual strategy. For a given strategy $S$, we denote by $S^d$ the dual strategy to $S$, where ``dual'' refers to switching the roles of W and B. Equivalently, we view a strategy $S$ as a pair of functions (one for each player) $f_S^1,f_S^2:P(X)\to X$ where $X$ is the set of hats being considered, and we view the two functions as taking the set of black hats on the partner's head as input and outputting the player's guess. Then, the dual strategy $S^d$ has as its functions $f_{S^d}^i(A) = f_S^i(X\setminus A)$.

\begin{remark}
In reference to the 4 optimal 2-player strategies, it is worth noting that $S_3 = S_1^d$, $S_2 = S_2^d$ (up to reordering of hats), while $V_2(p;S_4) = V_2(p;S_4^d)$, but we do not currently know whether $S_4^d=S_4$.
\end{remark}

The following lemma gives a formula for calculating the value of a dual strategy in terms of the value of the original strategy. For a given strategy $S$ for the hat game with probability $p$ of a B hat, let $\Prob_{S}(x_1,x_2)(p)$ be the probability that Player 1 chooses an $x_1$ hat on her head ($x_1 \in \{ W, B \}$) and Player 2 chooses an $x_2$ hat on his head.

\begin{lemma} \label{dual lemma}
Given any strategy $S$, let $S^d$ be its dual strategy. Then, for any $p \in (0,1)$,
\[
V_2(p; S^d) = 2p - 1 + V_2(1-p; S).
\]
\end{lemma}
\begin{proof}
First observe that by the definition of a dual strategy, $\Prob_{S^d}(B,B)(p) = \Prob_{S}(W,W)(q)$. This can be seen by pairing scenarios where all hat colors are reversed.

Observe that if the players utilize strategy $S$ when the probability of a black hat is $q$, then $\Prob_S(W,W)(q) + \Prob_S(W,B)(q) = 1-q = p$ is the probability of the first player selecting a white hat. Similarly, $\Prob_S(W,B)(q) + \Prob_S(B,B)(q) = q$ is the probability of the second player selecting a black hat.

Combining these observations, we have
\begin{align*}
V_2(p; S^d) & = \Prob_{S^d}(B,B)(p) \\
& = \Prob_{S}(W,W)(q) \\
& = p - \Prob_S(W,B)(q) \\
& = p - \left(q - \Prob_{S}(B,B)(q)\right) \\
& = p - q + \Prob_{S}(B,B)(q) \\
& = 2p - 1 + V_2(q; S).
\end{align*}
\end{proof}

\subsection{Proofs of Upper Bounds on $V_2(p)$}

In order to derive a general upper bound on $V_2(p)$ for any rational $p \in (0,1)$, we can apply the methods of Section~\ref{Sec:Matrix} to a $b \times b$ matrix where the first column is $a$ 1s followed by $b-a$ 0s, and the other columns are cyclic permutations of it. We 
compute upper bounds on the value of this matrix game,
which are in turn
upper bounds on the 
value $V_2(p)$ of the
hat game for $p=\frac{a}{b}$.

Below we restate Theorem~\ref{Thm:Upper} 
from the end of Section~\ref{Sec:Matrix} and then prove the theorem with the
help of a lemma.

\begin{reptheorem}{Thm:Upper}
For $p = \dfrac{a}{b}\leq \dfrac{1}{2}$, we have $\displaystyle V_2(p) \leq \frac{a}{b} - \left(1-\frac{a}{b}\right)^b \left(\frac{a}{b}\right)$.

For $p = \dfrac{a}{b}\geq \dfrac{1}{2}$, we have $\displaystyle V_2(p) \leq \frac{a}{b} - \left(1-\frac{a}{b}\right) \left(\frac{a}{b}\right)^b$.
\end{reptheorem}

Using
the notation from Section~\ref{Sec:n_players} 
(e.g., with $\Prob(BB)$ for the probability that each player picks a black hat),
we have
\[
\Prob(BB) = \Prob(B\ast) - \Prob(BW) .
\]

We know that $\Prob(B\ast) = p$. 
Thus all that remains is to bound the value of $\Prob(BW)$.

\begin{lemma}\label{Lem:Upper}
With the additional matrix-based information at the start of this section provided to the players, there is no strategy 
that wins with probability greater than
\[
\frac{a}{b} - \left(\frac{a}{b}\right)^b\left(1-\frac{a}{b}\right) .
\]
\end{lemma}

\begin{proof}
We use the observation above that
\[
\Prob(BB) = p - \Prob(BW).
\]
In particular, suppose that the players' agreed-upon strategy is that Player 1 will choose a hat from a fixed list $x_1,\ldots, x_k$ ($k\leq b$). It may happen that all of these hats are black, which happens with probability $p^k\geq p^b$, in which case Player 2 picks a white hat with probability $1-p$. This means that $\Prob(BW) \geq p^b(1-p)$. 
Therefore,
\[
V_2\left(\frac{a}{b}\right) \leq \frac{a}{b} - \left(\frac{a}{b}\right)^b\left(1-\frac{a}{b}\right).
\]
\end{proof}

To complete the proof of Theorem~\ref{Thm:Upper}, we observe that by 
the duality
discussed above, we also have
\[
V_2\left(\frac{a}{b}\right) \leq \frac{a}{b} - \left(1-\frac{a}{b}\right)^b\left(\frac{a}{b}\right).
\]
We compare the two bounds and discover that they are stronger on the intervals claimed in Theorem~\ref{Thm:Upper}.

\begin{remark}
The upper bound in Lemma~\ref{Lem:Upper} is sharp for the case when 
$p = \frac{b-1}{b}$
in the sense that one can describe a strategy for the 
$b \times b$
matrix game that succeeds with this probability. (However, this remains only an upper bound on $V_2(p)$ for the original hat game, in which the players do not receive hints.)
\end{remark}

\begin{remark}
Supposing that $V(p)$ is differentiable, using Theorem~\ref{Thm:Upper}, one can calculate an upper bound on the derivative of the function $V(p)$ at 0, and a lower bound at 1. Also, one can use the strategies described in Section~\ref{Sec:2p} to calculate a lower bound at 0 and upper at 1. One obtains $\frac{1}{2} \leq V'(0) \leq 1 - \frac{1}{e}$ and $1 + \frac{1}{e} \leq V'(1) \leq \frac{3}{2}$.
\end{remark}

The upper bounds on $V_2(p)$ in Theorem~\ref{Thm:Upper} are quite good for $p$ of the form 
$\frac{1}{b}$ or $\frac{b-1}{b}$,
as indicated in Figure~\ref{Fig:Bounds} 
(at the end of Section~\ref{Sec:Matrix})
by their proximity to the continuous lower-bound curve. However, they can be improved by using larger and less structured hint matrices. We give examples of this improvement below for the cases $p = 1/3$ and $p = 2/3$, whose upper bounds from Theorem~\ref{Thm:Upper} are $19/81 = 0.234567...$ and $46/81 = 0.567901...$, respectively. (The respective lower bounds from Theorem~\ref{Thm:Distinct} in Section~\ref{Sec:2p} are $0.205555...$ and $0.538888...$.)

The hint matrix
\[
U=
\begin{pmatrix}
1 & 0 & 0 & 1 & 0 & 0 & 0 & 0 & 0 & 0 & 1 & 1\\
1 & 0 & 0 & 0 & 1 & 0 & 0 & 0 & 1 & 1 & 0 & 0\\
1 & 0 & 0 & 0 & 0 & 1 & 1 & 1 & 0 & 0 & 0 & 0\\
0 & 1 & 0 & 1 & 0 & 0 & 0 & 1 & 0 & 1 & 0 & 0\\
0 & 1 & 0 & 0 & 1 & 0 & 1 & 0 & 0 & 0 & 0 & 1\\
0 & 1 & 0 & 0 & 0 & 1 & 0 & 0 & 1 & 0 & 1 & 0\\
0 & 0 & 1 & 1 & 0 & 0 & 1 & 0 & 1 & 0 & 0 & 0\\
0 & 0 & 1 & 0 & 1 & 0 & 0 & 1 & 0 & 0 & 1 & 0\\
0 & 0 & 1 & 0 & 0 & 1 & 0 & 0 & 0 & 1 & 0 & 1\\
\end{pmatrix}
\]
yields an upper bound of $0.221307...$ for $p=\frac{1}{3}$.

The hint matrix
\[
U=
\begin{pmatrix}
0 & 1 & 1 & 0 & 1 & 1 & 1 & 1 & 1 & 1 & 0 & 0\\
0 & 1 & 1 & 1 & 0 & 1 & 1 & 1 & 0 & 0 & 1 & 1\\
0 & 1 & 1 & 1 & 1 & 0 & 0 & 0 & 1 & 1 & 1 & 1\\
1 & 0 & 1 & 0 & 1 & 1 & 1 & 0 & 1 & 0 & 1 & 1\\
1 & 0 & 1 & 1 & 0 & 1 & 0 & 1 & 1 & 1 & 1 & 0\\
1 & 0 & 1 & 1 & 1 & 0 & 1 & 1 & 0 & 1 & 0 & 1\\
1 & 1 & 0 & 0 & 1 & 1 & 0 & 1 & 0 & 1 & 1 & 1\\
1 & 1 & 0 & 1 & 0 & 1 & 1 & 0 & 1 & 1 & 0 & 1\\
1 & 1 & 0 & 1 & 1 & 0 & 1 & 1 & 1 & 0 & 1 & 0\\
\end{pmatrix}
\]
yields an upper bound of $0.554641...$ for $p=\frac{2}{3}$.

\section{$n$-Player Results}\label{App:n_player_results}

\noindent
{\bf Note:} As in Section~\ref{Sec:n_players}, throughout Appendix~\ref{App:n_player_results}
we will think of the players as the wise men of Tanya Khovanova’s problem statement 
\cite{TK} and use singular masculine pronouns when referring to individual players.

\medskip

\subsection{Proof of Theorem~\ref{thm:ori}}\label{Subsec:ProofThmOri}

Below we prove the partial converse Theorem~\ref{thm:ori} from Section~\ref{Sec:n_players} to Peter Winkler's order-($1/ \log(n)$) strategy for $n$ players with arbitrary black-hat probability $p$.
As in Section~\ref{Sec:n_players}, we define $q := 1-p$ and $r := 1/q$.

\begin{reptheorem}{thm:ori}
Suppose that each player must choose the {\em lowest} level on his own head that has a black hat (with the team failing if any player has only white hats). Then, using $\tilde{V}_n(p)$ to refer to the optimal probability as $h \rightarrow \infty$ that the players succeed (for any given black-hat probability $p \in (0,1)$) under this more stringent requirement,  
for each $\epsilon$ with $0 < \epsilon \le 1/4$ we have
\[
\tilde{V}_n(p) \le (1 + \epsilon) / \log_r(n) \mbox{ for all sufficiently large } n.
\]
\end{reptheorem}

\begin{proof}
For $j = 1, \ldots, n$, let the random variable $Y_j$ be the lowest level on which Player $j$ has a black hat. We will restrict attention to the first $t = \lceil (1 - \epsilon / 2) \log_r(n) \rceil$ levels, and we will establish the desired upper bound on $\tilde{V}_n^{(h)}$ for all $h \ge t$, from which the bound will follow immediately for $\tilde{V}_n = \sup_h \tilde{V}_n^{(h)}$. (If any players have all-white $h$-hat stacks, the team would automatically lose the actual game, but such players will be considered to have $Y_j > t$ and will generously be exempted from having to guess at all for the purpose of this upper bound.)

An accomplice will uniformly at random choose a level $k$ from $\{1, \ldots, t\}$ and then uniformly at random choose a player ($\jstar$, say) from among those who happen to have $Y_j = k$ for the given realization. The accomplice will then inform the chosen player $\jstar$ of his special status but will not tell him his lowest-black-hat level $k$. Player $\jstar$ is the only player required to guess his level $Y_j$; the accomplice will tell the other $(n-1)$ players (including those with $Y_j > t$) their own values $Y_j$, and if any players have all-white $h$-hat stacks, they will be exempted from having to guess their values $Y_j$.

We argue
below that, with very high probability, all values 
$k \in \{1, \ldots, t\}$
occur and, in fact, occur about as many times as expected at random.
It then follows
that, even with this substantial help from the accomplice (which any or all of the players can ignore if they wish), 
the one player who is not told his own level still cannot 
pick the correct level with probability greater than $(1 + \epsilon) / \log_r(n)$ as $n \rightarrow \infty$.

Now we formalize the claim above. For each $k \in \{1, \ldots, t\}$, let $X_k$ be the number of players $j$ with $Y_j = k$. Then $\mu_k$, the expected value of $X_k$, satisfies
\[
\mu_k = n p q^{k-1} \ge n p q^{t-1} \ge p n^{\epsilon / 2} . 
\]
By a 2-sided multiplicative Chernoff bound for each value $k$, followed by a union bound over the $t$ possible values for $k$, we have, for each $\delta$ with $0 < \delta \le 1$,
\begin{align*}
\Prob( (1 - \delta) \mu_k \le X_k \le (1 + \delta) \mu_k \mbox{ for all } k \in \{1, \ldots, t\} ) & \ge 1 - 2t \exp(-(\delta^{2} / 3) p n^{\epsilon / 2}) \\
& = 1 - o(1/n^C) \mbox{ for every } C > 0 \\
& = 1 - o(1 / \log_r(n)) .
\end{align*}
Thus, for the purpose of proving Theorem~\ref{thm:ori}, we can neglect the probability that some value $k \in \{1, \ldots, t\}$ fails to occur or occurs with relative frequency significantly different from its expected value.

Now the chosen player $\jstar$ sees the lowest-black-hat level $Y_j$ for each of the other players $j$ and knows how he was chosen; this allows him to compute a well-defined posterior probability distribution for his own level $Y_{\jstar} \in \{1, \ldots, t\}$. Since $X_k / \mu_k \in [1 - \delta, 1 + \delta]$ for all $k \in \{1, \ldots, t\}$ with all but asymptotically negligible probability, it follows readily from Bayes' theorem that the ratios of posterior probabilities
\[
\Prob_{\mbox{post}} (Y_{\jstar} = k_1) / \Prob_{\mbox{post}} (Y_{\jstar} = k_2)
\]
are in $[(1 - \delta)/(1 + \delta), (1 + \delta)/(1 - \delta)]$ for all $k_1, k_2 \in \{1, \ldots, t\}$. If we let $\delta = \epsilon / 16$, say, it then follows readily that, for each $k \in \{1, \ldots, t\}$,
\[
\Prob_{\mbox{post}}(Y_{\jstar} = k) \le \frac{1 + \epsilon /4}{t} \le \frac{(1 + \frac{\epsilon}{4}) / (1 - \frac{\epsilon}{2})}{\log_r(n)} .
\]
Recalling that $0 < \epsilon \le 1/4$ and taking into account the asymptotically negligible probability of atypical events, we find that, even with the help from the accomplice (which cannot hurt the players, since they are free to ignore extra information),
\[
\tilde{V}_n(p) \le \frac{1 + \epsilon}{\log_r n}
\]
for all sufficiently large $n$.

\end{proof}

\subsection{Proof That Mis{\`e}re Success Probability Is Bounded away from 1} \label{app:subsec:W_n_lt_1}

Now we complete the proof of Theorem~\ref{thm:misere} from Section~\ref{Sec:n_players} to show that $W_n$, the probability of success for $n$ players in the mis{\`e}re game, is bounded away from 1 for each $n$. We begin by recalling the definition of the pairs of integers $(r_j, s_j)$ for $j \ge 1$ as below; these will be the dimensions of ``hint matrices" given to the various players for the $n$-player mis{\`e}re game.

\begin{definition}
Let $(r_1, s_1) = (2,2), \ (r_2, s_2) = (4, 6)$, and for each $k \ge 3$, define 
\[
r_k = 2 \prod_{j=1}^{k-1} s_j \mbox{ and } s_k = \binom{r_k}{r_k / 2}.
\]
\end{definition}

\begin{reptheorem}{thm:misere}
Letting $W_n = \sup_S \Prob(\mbox{At least 1 of $n$ players chooses a white hat})$ over all possible $n$-player strategies as $h$, the number of hats per head, grows to infinity, we have the following:
\[
W_2 = 1/2 + V_2 \in [17/20, 193/224] \mbox{ (hence in [0.85, 0.8616...])}
\]
and, for $n \ge 3$,
\[
1 - (1/2)^n \le W_n \le 1 - \frac{1}{(r_{n} / 2) 2^{(r_{n} / 2)} }.
\]
\end{reptheorem}

\begin{proof}
The bounds on $W_2$ and the lower bounds on $W_n$ for all $n \ge 3$ were already established in Section~\ref{Sec:n_players}. The upper bounds follow from a ``hint" technique generalizing the upper-bounding techniques used earlier for $V_2$ and $V_3$. 

For $1 \le j \le n-1$, Player $j$ is given an $r_j \times s_j$ hint matrix with $r_j$ and $s_j$ as defined just before the statement of the current theorem. Player $n$ is given no hint, but by the usual argument, he can 
without loss of generality
restrict attention to (at most) the first $r_n / 2$ levels on his own head, since this is the product of the number of columns in the other $n-1$ players' hint matrices, which is the total number of distinguishable situations in which Player $n$ can find himself. With probability at least $1/(2^{r_n / 2})$ (strictly greater than this if Player $n$ does not actually use all possible $r_n / 2$ levels on his own head), he will have black hats on all of the levels from which he chooses, and Players 1 through $n-1$ will all know when they are in this situation. In this case, it is up to the first $n-1$ players to choose at least 1 white hat. 

Now, proceeding inductively downstream from Player $n-1$ through Player 2, conditioning on what is seen on the heads of all the upstream players $k+1, \ldots, n$, each Player $k$ sees one of $s_1 \cdot s_2 \ldots \cdot s_{k-1} = r_k / 2$ possible joint column choices for Players 1 through $k-1$. Even if Player $k$ assigns a different row of his own hint matrix to each of these $r_k / 2$ distinguishable downstream possibilities, one of his $\binom{r_k}{r_k / 2}$ columns will contain 
1s in all $r_k / 2$ of these rows. (If Player $k$ sometimes assigns the same row to different distinguishable downstream observations, there will be multiple columns of his hint matrix that contain 
1s 
in all rows that he actually uses.) 

Thus, conditioned on whatever Player $k$ observes upstream and downstream (and whatever strategy he has committed himself to), with probability at least $1 / s_k$ he will have been assigned a column by the referee that forces him to choose a black hat, inductively leaving the downstream Players 1 through $k-1$ with the responsibility of choosing at least one white hat. Finally, if Players 2 through $n$ have all been assigned these most unfavorable columns by the referee, Player 1 will know this fact and will have the burden of choosing a white hat on his own head. However, his 2 possible hint columns are equally probable and differ from each other on every level, so Player 1 will fail with probability 1/2. Multiplying all $n$ of the players' respective conditional failure probabilities together, we see that they must lose the mis{\`e}re game with probability at least
\[
(1/2^{r_n / 2}) \cdot (1/s_{n-1}) (1/s_{n-2}) \ldots (1/s_1).
\]

Since $r_n$ is defined as $2 s_1 s_2 \ldots s_{n-1}$, the claimed upper bound on $W_n$ follows immediately. 
\end{proof}

\subsection{Proof that Mis{\`e}re Success Probability $W_n$ Decreases in $n$}\label{app:subsec:W_n_decreases}

The next result, mentioned in Section~\ref{Sec:n_players}, shows that the optimal success probability $W_n$ for the $n$-player mis{\`e}re game with $p = 1/2$ (in which {\em at least} one player must point to a {\em white} hat) is strictly increasing in $n$.

\begin{theorem}\label{thm:misere_increasing}
The optimal mis{\`e}re success probabilities satisfy
$W_{n+1} > W_n$ for all $n \ge 1$, with $\lim_{n \rightarrow \infty} W_n = 1$.
\end{theorem}

\begin{proof}
The fact that $\lim_{n \rightarrow \infty} W_n = 1$ follows immediately from Theorem~\ref{thm:misere}. Now, much as in the proof of the previous theorem, we have
\[
\Prob(B \ldots B,*) = \Prob(B \ldots B,B) + \Prob(B \ldots B,W)
\]
for any $(n+1)$-player strategy, so
\[
1 - \Prob(B \ldots B,*) = (1 - \Prob(B \ldots B,B)) + (1 - \Prob(B \ldots B,W)) - 1
\]
for any $(n+1)$-player strategy. Thus
\[
\sup \ \{ 1-\Prob(B \ldots B,*) \} = \sup \ \{(1-\Prob(B \ldots B,B)) + (1-\Prob(B \ldots B,W)) \} - 1,
\]
where the supremum is over all $(n+1)$-player strategies. Then
\[
\sup \ \{ 1 - \Prob(B \ldots B,*) \} \le \sup \ \{ 1 - \Prob(B \ldots B,B) \} + \sup \ \{ 1 - \Prob(B \ldots B,W) \} - 1.
\] 
Since the maximal probability of avoiding $(B \ldots B,B)$ is the same as the maximal probability of avoiding $(B \ldots B,W)$ when $p = 1/2$, we have
\[
W_n \le 2 W_{n+1} - 1,
\]
so $W_{n+1} \ge (1 + W_n) / 2 > (W_n + W_n)/2$, where the last inequality is because $W_n < 1$ for each $n$.

Thus we have the strict inequality $W_{n+1} > W_n$ for all $n \ge 1$. 
\end{proof}

\subsection{Details of Best Strategies Known for $n$ Players}\label{App:subsec:details_best_strategies}

Now we give details about the best lower bounds known for $V_3$ through $V_{12}$ (all for black-hat probability $p = 1/2$) that were omitted in Section~\ref{Sec:best_lower_bounds}. We discuss strategies found by hill climbing for $V_3$ and $V_4$ and describe generalizations of the basic Winkler strategy from Theorem~\ref{thm:winkler} in Section~\ref{Sec:n_players}. Finally, we analyze (recursive) white-reset and black-reset enhancements for $t$-hat tiers that improve all of the strategies above.

In an email sent in 2014 to various hats enthusiasts, Jay-C Reyes and Larry Carter reported what were then the best lower bounds known on $V_3$ and $V_4$ for $p = 1/2$. Their bounds were constructive and came from hill-climbing on symmetric strategies for 4 hats per player. (Since we will later extend these strategies by considering 4-hat {\em tiers}, 
we use
$t$ rather than $h$ to refer to the number of hats per player in each basic strategy.) The Carter-Reyes 4-player result leads to what is still the best known lower bound for $V_4$. Their 3-player search has since been adapted to consider $t = 3, 4, 5,$ or $6$ hats per player. Before incorporating the ``reset" enhancements, we obtain the following lower bounds from these searches:
\[
\begin{array}{ccccccc}
V_3 & \ge & 9120/(2^5)^3 & = & 9120/32768 & = & 0.278320... ,\\
V_4 & \ge & 14845/(2^4)^4 & = & 14845/65536 & = & 0.226516... .
\end{array}
\]

\begin{remark} 
The 3-player strategy, which uses 5 levels, can be described by a symmetric $32 \times 32$ matrix with values in $\{1,2,3,4,5\}$ that is used by all 3 players. The 4-player strategy, which uses 4 levels, can be described by a $16 \times 16 \times 16$ symmetric tensor with values in $\{1,2,3,4\}$ that is used by all 4 players.
\end{remark}

For $n \ge 5$, our best lower bounds on $V_n$ come from Winkler-style strategies in which players focus on the first $t \approx \log_2(n)$ levels. With $Y_j$ defined as the lowest level on which Player $j$ has a black hat, players using the original Winkler strategy hoped that $Y_1 + \ldots + Y_n$ would have some particular residue (mod $t$) with probability as much above the guaranteed $1/t$ as possible. This best-residue probability can be improved a little if each player is allowed to apply some permutation $\pi_j$ to his value $Y_j$ before the values are summed. (The players' permutations on $\{1, \ldots, t\}$ can be different for different players, as long as they are fixed during the strategy session.) For example, for the usual case $p = 1/2$, one appears to do better by computing the alternating sum $(Y_1 - Y_2 + Y_3 - \ldots + (-1)^{n+1}Y_n) \mod t$ than the straight sum (mod $t$), since the alternating sum concentrates the probability mass more effectively.

Still more generally, one can use $t \times t$ Latin squares to ``add in" one player's level at a time to the ``running sum," possibly permuting the output symbols after each new player's level is folded in. When $t = 4$, for example, one does better by using Latin squares corresponding to the Klein 4-group than by using squares corresponding to addition or subtraction (mod 4).

The bounds for $5 \le n \le 8$ all come from applying the Winkler idea to 4-hat tiers, representing each level within a tier by a dibit in $\{00, 01, 10, 11\}$, and XORing the $n$ dibits corresponding to the lowest level in the tier (if any) on which each player has a black hat, resetting to the next tier if necessary. The players hope that they will each have at least one black hat within the 4-hat tier and that the XOR of the $n$ resulting dibits will be $00$, and they each choose the corresponding one of 4 levels on their own heads. It turns out that for this strategy, the players win if and only the mod-2 sum of all their dibits really is $00$. 

However, for many strategies (most notably, for the original strategy of straight addition of lowest-black-hat levels (mod $t$)), there is ``secondary success probability'' (or ``bonus'' probability) coming from cases in which the actual ``sum'' of the $n$ lowest levels differs from the targeted value but the players nonetheless serendipitously each point to a black hat. (In fact, because of this bonus probability, it turns out that straight mod-$t$ addition almost always yields more overall success than alternating addition and subtraction (mod $t$), even though the latter almost always yields higher ``primary'' success probability than the former.) We have accounted for this bonus probability in all of our best known lower bounds in Section~\ref{Sec:best_lower_bounds}.

The XOR strategies with $t = 4$ equivalently use the Latin square
\[
\begin{pmatrix}
0 & 1 & 2 & 3\\
1 & 0 & 3 & 2 \\
2 & 3 & 0 & 1 \\
3 & 2 & 1 & 0
\end{pmatrix}
\]
associated with $Z_2 \times Z_2$ to combine values corresponding to the lowest-black-hat level within the tier for each successive player. These strategies outperform strategies based on mod-4 addition and subtraction (and many permutation-based generalizations thereof), apparently because the Latin square associated with $Z_2 \times Z_2$ concentrates probability within the first two categories more effectively than the Latin squares associated with those other arithmetic operations.

For $9 \le n \le 12$, our best strategies use 5-hat tiers (i.e., $t = 5$). There are two distinct isotopy classes of $5 \times 5$ Latin squares (as there are for $4 \times 4$ Latin squares), and once again, it appears that the isotopy class {\em not} associated with mod-$t$ addition or subtraction does a better job of concentrating probability within 2 of the $t$ categories. The $5 \times 5$ isotopy class that we found to work best corresponds to a nonassociative quasigroup (in fact, to a loop). Our best results are probably not optimal even within the class of strategies we considered, since we used a greedy search algorithm, but they all begin by combining the first two mod-5 values using the Latin square
\[
\begin{pmatrix}
0 & 1 & 2 & 3 & 4\\
1 & 0 & 3 & 4 & 2 \\
2 & 4 & 0 & 1 & 3 \\
3 & 2 & 4 & 0 & 1 \\
4 & 3 & 1 & 2 & 0
\end{pmatrix} .
\]

Later mod-5 values are folded in using Latin squares isotopic to this first square, but not the same square (or quasigroup) for each new player.

All of the strategies above are improved slightly by working with $t$-hat ``tiers'' and recursively ``resetting'' (shifting up $t$ levels at a time) whenever a player sees another player with an all-white stack of hats within the current tier. In any such situation, the players would certainly lose without resetting, and if there is only one player with an all-white stack within the current tier, he will fail to get the memo, so the team will still lose. However, if 2 or more players have all-white stacks in the current $t$-level tier, all $n$ players will reset together and give themselves an independent chance of winning at the next tier. By recursively resetting, the players succeed for all placements of hats on the lowest $t$ levels that would have won without resetting, but now the denominator of their success probability is reduced from $2^{tn}$ to
\[
(2^t - 1)^n + \binom{n}{1} (2^t - 1)^{n-1} .
\]

One can refine the reset strategy a bit. If any player (call him Player $j$) sees only {\em black} hats within the first $t$-level tier on {\em all} of the other $n-1$ heads, he can reset (shifting up $t$ levels on all of the stacks that he sees).

If he has neither all-black nor all-white on his own first $t$ levels, then the other $n-1$ players will stay on the first $t$ levels and will all be guaranteed to guess correctly. In this case, the team wins with average conditional probability 1/2 whether Player $j$ (the non-monochromatic player) guesses from his first $t$-level tier or from any other $t$-level tier. 

If Player $j$ has only white hats in his first tier, then the other $n-1$ players will also reset (since they see his all-white stack), and in this case, the team will get a fresh start at the next tier of $t$ levels, whereas they would have lost for sure if they had used no resets or only the reset-on-white strategy. If Player $j$ has all black hats in his first $t$-level tier, then all $n$ players have all-black first tiers, so they will all reset to the next tier and give themselves a fresh start of winning with conditional probability $V_n$ rather than the conditional probability of 1 that they would have enjoyed if they had stayed put.

We see that this reset-on-black strategy will change the probability of winning (with respect to the earlier strategy of resetting only when at least one player has an all-white first tier) only when everyone has a monochromatic first tier, with at most one player having an all-white first tier. This situation occurs with probability $(n+1) (1/2)^{nt}$, and the conditional net gain in success probability in this case is at least
\[
\frac{n}{n+1} (\shortunderline{V}_n - 0) + \frac{1}{n+1} (\shortunderline{V}_n - 1) = \shortunderline{V}_n - 1/(n+1),
\]
where $\shortunderline{V}_n$
is any valid lower bound on $V_n$. Thus as long as we have a starting strategy that achieves success probability strictly greater than $1/(n+1)$, as we do for all $n \ge 2$, this augmented resetting strategy helps. 

When this reset-on-black policy is implemented recursively, the denominator of the players' success probability is reduced by $n+1$, to
\[
(2^t - 1)^n + n (2^t - 1)^{n-1} - (n+1) ,
\]
and the numerator is reduced by 1.

With resetting incorporated into the basic $t$-hat strategies, we obtain the following lower bounds on $V_3$ and $V_4$.

\begin{lemma}
For the 3-player and 4-player games, we have
\[
 V_3 \ge 9119/32670 = 0.279124... ~~~\mbox{and}~~~~ V_4 \ge 14844/64120 = 0.231503....
\]
\end{lemma}

\begin{proof} 
For $n = 3$, we used 5 levels and found a symmetric strategy (described by a symmetric $32 \times 32$ matrix with values in $\{1, 2, 3, 4, 5\}$ that is used by all 3 players) that wins 9120 times out of 32768. With white and black resets, this yields success probability
\[
\frac{9120-1}{31^3 + \binom{3}{1} (31)^2 - (3+1)} = \frac{9119}{(31)^2 (31 +3) -4} = \frac{9119}{32670} .
\]

For $n = 4$, Reyes and Carter used 4 levels and found a symmetric strategy that wins with probability $14845/65536$, which is improved to $14844/64120$ with recursive resetting. 
\end{proof}

For large values of $n$, the best lower bounds we know how to achieve are only very slightly above the obvious lower bounds from the basic Winkler strategy with $t$-hat tiers and a reset to the next tier if a player sees at least one other player who has only white hats in his first tier. (Resetting when one sees only {\em black} hats on {\em all} other players' heads provides much less help for large $n$.) With just the reset on white, we obtain
\[
V_n \ge \max_{t} \frac{1}{(2^t - 1)^n + n (2^t - 1)^{n-1}} \ \frac{(2^t - 1)^n}{t} = \max_t \frac{2^t - 1}{2^t - 1 + n} \ \frac{1}{t},
\]
which decreases asymptotically as $1 / \log_2(n)$, essentially as argued in Theorem~\ref{thm:winkler}.

\end{document}